%% file: bifurcations.tex
\documentclass[12pt,twoside]{mymathart}
\usepackage[theorems]{mymathmacros}

\usepackage{ifthen}
\usepackage{graphicx,subfigure}
\usepackage{hyperref}

\usepackage[english,german]{babel}

\input{symboldefs.tex}

\newcommand{\IR}{\Gamma_{W,\hspace{.5pt}h}}
\newcommand{\IRH}[1]{\Gamma_{W,\hspace{.5pt}#1}}
\newcommand{\IRW}[1]{\Gamma_{#1,\hspace{.5pt}h}}
\newcommand{\IRWH}[2]{\Gamma_{#1,\hspace{.5pt}#2}}
\newcommand{\PR}{G_{\s}}
\newcommand{\PRS}{G}
\newcommand{\M}{\mathcal{M}}

\newcommand{\PBM}{\chi_n}

\begin{document}
\selectlanguage{english}
\title{Bifurcations in the Space of Exponential Maps}

\author{Lasse Rempe}
\address{Department of Mathematical Sciences, University of Liverpool, L69 7ZL,
United Kingdom}
\email{l.rempe@liverpool.ac.uk}
\thanks{The first author was supported in part
 by a fellowship of the 
 German Academic Exchange Service and
 by the German-Israeli Foundation
 for Scientific Research and Development,
 grant no.\ G-643-117.6/1999}

\author{Dierk Schleicher}
\address{International University Bremen, P.O.~Box 750 561, 
28725 Bremen, Germany}
\email{dierk@iu-bremen.de}

\subjclass{Primary 37F10; Secondary 30D05}

\date{\today}


\begin{abstract}
 This article investigates the
  parameter space of the exponential family $z\mapsto
 \exp(z)+\kappa$. We prove that the boundary (in $\C$) of every hyperbolic
 component is a Jordan arc, as conjectured by Eremenko and Lyubich as
 well as Baker
 and Rippon. In fact, we prove the stronger statement that the
 exponential bifurcation locus is connected in $\C$, which is an analog
 of Douady and Hubbard's
  celebrated theorem that the Mandelbrot set is connected.
  We show furthermore 
  that $\infty$ is not accessible through any nonhyperbolic
  (``queer'') stable component. The main part of the argument consists
  of demonstrating a general
  ``Squeezing Lemma'', which controls the structure of parameter space
  near infinity. We also prove a second conjecture of Eremenko and
  Lyubich concerning bifurcation trees of hyperbolic components.
\end{abstract}

\maketitle

\tableofcontents


\section{Introduction}

This article is one step in a program to describe the dynamics of the family
 of
 exponential maps $\Ek:z\mapsto \exp(z)+\kappa$ and the structure of its
 parameter space. This simplest family among transcendental entire maps
 has
 received special attention over the years, much like
 the quadratic family has among polynomials.
 A good 
 understanding of exponential dynamics has often 
 entailed progress in the study of
 more general classes of entire functions. 

Such study
 is important not only in its own right, but also because 
 the iteration
 of transcendental functions has links to many other areas in 
  dynamical systems and
 function theory. We content ourselves with giving some examples. 
 Features of exponential dynamics appear in the study of parabolic
  implosion \cite{mitsudim}; similarly, recent results on the hyperbolic
  dimension of Feigenbaum-type Julia sets \cite{avilalyubichfeigenbaum2}
  resemble 
  those obtained for hyperbolic exponential maps \cite{urbanskizdunik1}. 
  The family
  $\lambda te^{-t}$, a close relative of the exponential, is
  the second simplest model in population dynamics (the first being the
  logistic family), and the standard Arnol'd family of circle
  maps (a family of self-maps of
  the punctured plane), plays a prominent role in the renormalization of critical
  circle maps. Finally, anyone interested in finding
  the roots of an entire function should consider studying its
  Newton method, i.e., iteration of a transcendental meromorphic
  function.

 The reason that the exponential family is a natural candidate to begin with
  is the same that has made the quadratic family a favorite object of study:
  in both cases the maps possess only  one \emph{singular value}.
  The singular values (i.e., the critical
  and asymptotic values) of a function play an important role in the study of
  its dynamics: a restriction on the number of
  singular values generally limits the amount of different dynamical features
  that can appear for the same map. Therefore, the simplest non-trivial
  parameter space of holomorphic functions is given by the quadratic family,
  in which each function has only a single simple critical
  point in $\C$. Similarly,  
  exponential maps are the only transcendental maps with only
  one singular value in $\C$ (see, e.g., 
  \cite[Appendix D]{jacktwocriticalpoints}).
  Furthermore, the exponential family is the limit of the families of
  unicritical polynomials, parametrized as
  $z\mapsto (1+\frac{z}{d})^d+c$. This makes it an excellent candidate
  for applying methods that have proved useful in the study of these
  polynomials, as first developed for the Mandelbrot set $\M$
  in the
  famous Orsay Notes \cite{orsay}. 

 The first major result concerning the Mandelbrot set which Douady and
  Hubbard proved in these notes
  \cite[Th\'eor\`eme 1 in Section 3 of Expos\'e VIII]{orsay} 
  states that $\M$, or equivalently the 
  bifurcation locus in the
  space of quadratic polynomials, is
  connected. We prove the analogous fact for exponential parameter space.
  (See Section \ref{sec:boundary} for definitions.) 

 \begin{thm}[Bifurcation Locus is Connected]
    \label{mainthm:bifurcationlocus}
  The bifurcation locus  $\B$
   in the space of exponential maps is a connected subset of $\C$.
 \end{thm}

 To prove this, we need to show that no \emph{stable component}
  (that is, a component of $\C\setminus\B$) 
  can disconnect
  the parameter plane $\C$. 
  As in other complex analytic parameter spaces (such as the Mandelbrot set), 
  the
  most prominent open problem concerning exponential maps
   is to show that every such
   component is a
  \emph{hyperbolic component} (e.g., a component of the open subset
  consisting of maps which have an attracting periodic orbit.) 
  Therefore, as in the polynomial cases,
  a special interest lies
  in studying this type of component.  
  It is well-known \cite[Theorem 12]{alexmisha} that
  every hyperbolic component $W$ is simply connected and
  the multiplier map $\mu:W\to\Ds$, which maps each parameter to the 
  multiplier of its unique attracting cycle, is a universal
  covering. Therefore, there is a conformal isomorphism $\Phi_W\colon
  W\to\H$ with $\mu=\exp\circ\Phi_W$, where $\H$ is the left
  half plane.
  $\Phi_W$ is uniquely
  defined by this condition up to addition of $2\pi i\Z$; a preferred
  choice for $\Phi_W$ was described in \cite{habil,expattracting};
  compare Proposition \ref{prop:bifurcationstructure}.

\begin{figure}
 \center
 \includegraphics[height=0.87\textheight]{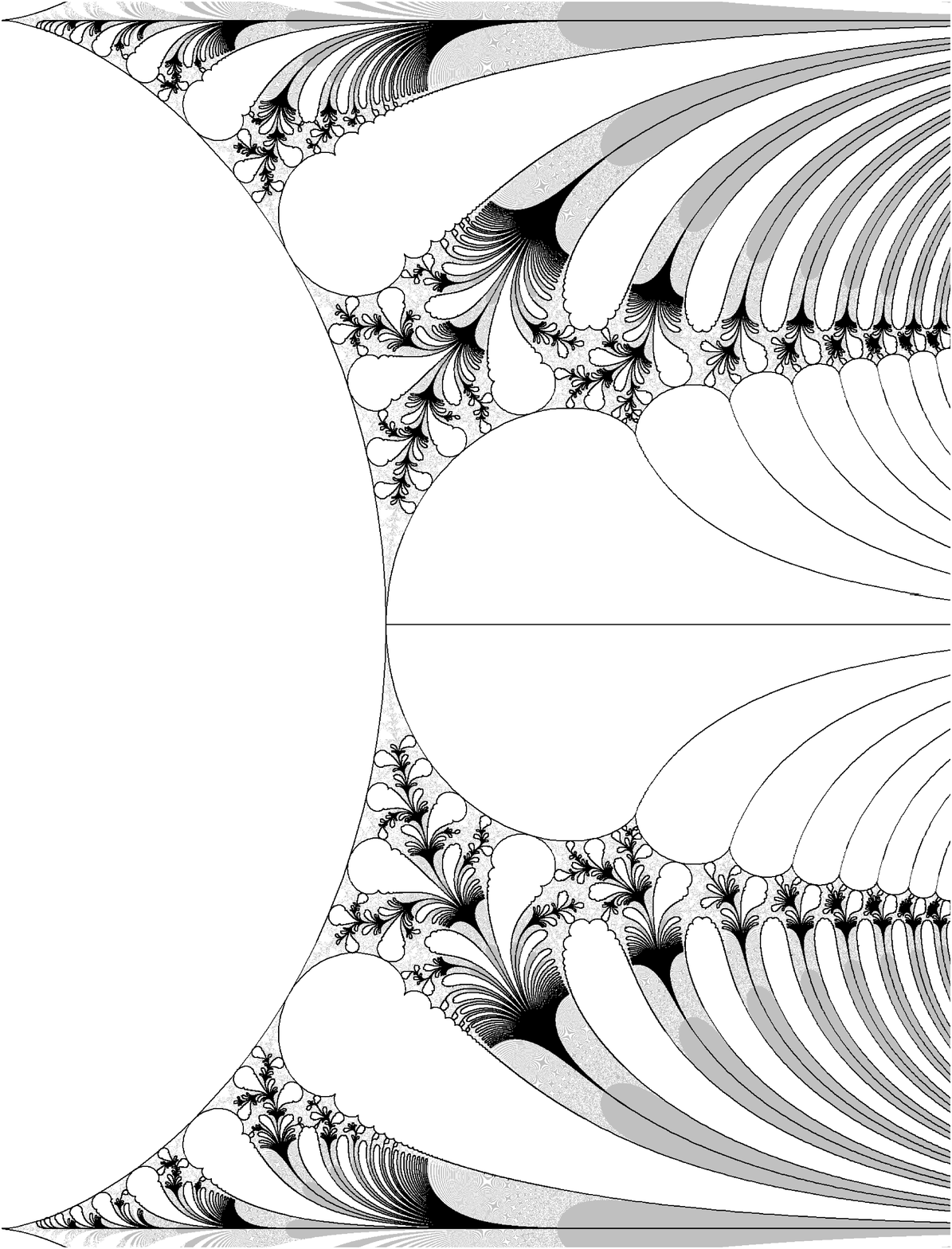}
 \caption[Hyperbolic Components]{Several hyperbolic components in the
 strip $\{\im \kappa \in [0,2\pi]\}$. Within the period two component
 in the center of the picture, the curves corresponding to real
 positive multipliers, i.e.~the components of the set
 $\mu^{-1}\bigl((0,1)\bigr)$, have been drawn in.\label{fig:hyperboliccomponents}}
\end{figure}

In the 1980s, both Baker and Rippon \cite{bakerexp} and 
 Eremenko and Lyubich \cite{alexmisharussian}
 independently conjectured
 that the boundary of a hyperbolic component is always a
 Jordan curve through infinity (see Figure
 \ref{fig:hyperboliccomponents}). 
 (Here and throughout the article, all closures
 are understood to be taken on the Riemann sphere $\Ch$.)
 The following theorem
 confirms this conjecture, and hence shows that a hyperbolic component cannot 
 disconnect the exponential bifurcation locus. 

\begin{thm}[Boundaries of Hyperbolic Components] \label{mainthm:boundary}
 Let $W$ be a hyperbolic component, and let $\Phi_W:W\to\H$ be a
 conformal isomorphism with $\mu=\exp\circ\Phi_W$. Then $W$ extends to
 a homeomorphism $\Phi_W:\cl{W}\to \cl{\H}$ with
 $\Phi_W(\infty)=\infty$. In particular, $\partial W\cap\C$ is a
 Jordan arc tending to $\infty$ in both directions.
\end{thm}

It is well-known (see Section \ref{sec:boundary}) that the map 
  $\Psi_W := \Phi_W^{-1}$ extends to
  a continuous map $\Psi_W:\cl{\H}\to\cl{W}$ with
  $\Psi_W(\infty)=\infty$. Thus the main difficulty in Theorem
  \ref{mainthm:boundary} lies in
  showing that any curve $\gamma:[0,\infty)\to W$
  with $|\gamma(t)|\to\infty$ satisfies
  $\Phi_W(\gamma(t))\to\infty$, or equivalently lies in the same
  homotopy class of $W$ as $\Psi_W\bigl((-\infty,0)\bigr)$. 
  The basic idea is to use the
  structure of parameter space near $\infty$ to exclude any other
  direction in which $\gamma$ could tend to $\infty$. This is the same
  technique used in the original (unpublished) proof of Theorem
  \ref{mainthm:boundary} \cite{habil}, as outlined in \cite{cras}. 
  (This proof, however, used
  landing properties of periodic parameter rays which were also
  established in \cite{habil}; our proof does not require any results
  of this sort.) In fact, we extend Theorem \ref{mainthm:boundary} to the
  following remarkable result, which shows that there are only two
  possibilities for a curve to infinity in
  exponential parameter space which does not intersect boundaries of
  hyperbolic components.

\begin{thm}[Squeezing Lemma] \label{thm:squeezing1}
 Let $\gamma:[0,\infty)\to\C$ be a curve in parameter space with
  $\lim_{t\to\infty}|\gamma(t)|=\infty$. Suppose that, for all $t$,
  the map
  $E_{\gamma(t)}$ does not have an indifferent periodic orbit. Then either
  \begin{itemize}
   \item for all sufficiently large $t$, 
          the singular value of $E_{\gamma(t)}$ escapes
          to $\infty$ under iteration, or
   \item $\gamma$ lies in a hyperbolic component $W$ and
          $\lim_{\t\to\infty}\Phi_W(\gamma(t))\to \infty$.
  \end{itemize}
\end{thm}

 The Squeezing Lemma not only proves Theorem \ref{mainthm:boundary}, but also
  implies Theorem \ref{mainthm:bifurcationlocus} on the connectedness
  of the bifurcation locus. Indeed, a (hypothetical) non-hyperbolic 
  stable component cannot contain any attracting, indifferent or
  escaping parameters. Hence we can apply the Squeezing Lemma to such 
  a component.

 \begin{cor}[Nonaccesibility of $\infty$ in Queer Components]
  \label{cor:queer}
  Suppose that $W$ is a non-hyperbolic stable component. If $W$ is
  unbounded, then $\infty$ is not accessible from $W$. 
  In particular, $\C\setminus W$ is connected. 
 \end{cor}
 
 The Squeezing Lemma is useful beyond its applications in this
  article. For example,
  a more precise version (Theorem \ref{thm:squeezinglemma})
  can be used to show that parameter rays
  cannot land at $\infty$. This is important in order to complete the
  classification of escaping parameters begun in \cite{markus}; see
  \cite{markuslassedierk}. 
  It does not seem unreasonable to expect that a stronger version of the 
  Squeezing Lemma is true 
  which could be used to prove
  that any non-hyperbolic component must be bounded. This
  is the
  subject of ongoing research.

\subsection*{Historical remarks}
 Parameter space for exponential maps was first
  studied in the articles \cite{bakerexp} and
  \cite{alexmisharussian} 
  (the latter
  treated general parameter spaces of transcendental entire
  maps which have finitely
  many singular values), and \cite{devaneybifurcation}.
  The preprint \cite{dgh} (recently published as \cite{dghnew1,dghnew2})
  was the first to
  view the exponential family as a limit of the unicritical polynomial
  families. 

 Since these seminal papers, exponential maps have been studied from many
  different points of view, and by many 
  authors; giving justice to all of them is beyond the scope of this 
  introduction. Among recent results about exponential parameter space,
  the most relevant to our current investigation is the
  existence and classification of hyperbolic exponential maps
  \cite{expattracting,dfj}; 
  compare Theorem \ref{thm:hyperboliccomponents} below. 

  We should also mention that there are a number of
   recent theorems (e.g.\ 
   \cite{nonrecurrentmero,urbanskizdunik4}) 
   which imply that exponential maps with certain dynamical 
   properties are contained in the bifurcation locus. These 
   interesting results 
   prove that non-hyperbolic components cannot be situated at such 
   parameters, 
   but do not make any statements regarding 
   the shape of arbitrary stable components. Hence they do not impact 
   directly on our Theorems \ref{mainthm:bifurcationlocus} and 
   \ref{mainthm:boundary}. 

\subsection*{Structure of the article}
 The structure of this article is as follows. We start by using the
  Squeezing Lemma to prove Theorems \ref{mainthm:bifurcationlocus} 
  and \ref{mainthm:boundary} in Section
  \ref{sec:boundary}. In Section \ref{sec:combinatorics}, we review
  results from \cite{expescaping} and \cite{expattracting}, and
  formulate a more precise version of the Squeezing Lemma (Theorem
  \ref{thm:squeezinglemma}). Its statement breaks down into three
  parts, which are proved in the rest of the article.
  The first of these can be
  proved by an elementary argument (Section
  \ref{sec:squeezingb}). The other two parts (proved in Section
  \ref{sec:squeezingc} and \ref{sec:squeezinga}) require information on the
  structure of \emph{wakes} in exponential parameter space. These
  results from \cite{expcombinatorics} are reviewed in Section
  \ref{sec:wakes}.

 Finally, in Section \ref{sec:furtherresults}, we prove a second
 conjecture of Eremenko and Lyubich from \cite{alexmisharussian}:
 there are infinitely many \emph{bifurcation trees} of hyperbolic
 components. We also show that the boundary of a hyperbolic component
 $W$ is an analytic curve, with the possible exception of cusps in the
 points of $\Phi_W^{-1}(\Z)$. Furthermore, we 
 indicate further
 developments which are related to our results, such as the landing of
 periodic parameter (and dynamic) rays.

\subsection*{Some remarks on notation}
 We have chosen to parametrize our
  exponential maps as $z\mapsto \Ek(z)=\exp(z)+\kappa$. Traditionally,
  they have often been parametrized as $\lambda\exp$, which is
  conjugate to $\Ek$ if $\lambda=\exp(\kappa)$. We prefer our
  parametrization mainly because the behavior of exponential maps at
  $\infty$, and in particular the asymptotics of dynamic rays, do not
  depend on the parameter in this parametrization. Note that this is
  also the case in the usual parametrization of quadratic polynomials
  as $z\mapsto z^2+c$.
  Another conceptual advantage is that the parameter $\kappa$ equals the
  singular value of $\Ek$, and thus the picture
  in the parameter plane reflects the situation in the dynamical
  plane. Finally, our parametrization causes our statements on
  hyperbolic components to have fewer exceptions for the components of
  periods $1$ and $2$.
  Note that $\Ek$ and
  $E_{\kappa'}$ are conformally conjugate if and
  only if $\kappa-\kappa'\in 2\pi i \Z$.

 The
 function $F:(0,\infty)\to(0,\infty), t\mapsto e^t-1$ 
 will be fixed throughout the article as a
 model function for exponential growth.
 If $\gamma:[0,\infty)\to\C$ is a curve, we shall say that
  $\lim_{t\to\infty} \gamma(t)=+\infty$ if $\re
 \gamma(t)\to +\infty$ and $\im \gamma$ is bounded; analogously for
 $-\infty$. The $n$-th iterate of any function $f$ will be denoted by
 $f^n$. As already mentioned, the closure and boundary of a subset of
 the plane will always be taken on the Riemann sphere
 $\Ch$. We
 conclude any
 proof and any result which immediately follows from previously proved
 theorems by the symbol $\blacksquare$\hspace{2pt}. A result which is cited without
 proof is concluded by $\square$.

\subsection*{Acknowledgements}
We would like to thank Walter Bergweiler, Alexandre Eremenko,
 Mikhail Lyubich, John Milnor, Rodrigo Perez, Phil Rippon,
 Juan Rivera-Letelier and especially Markus F\"orster for many helpful
 discussions and comments, 
 and the Institute of Mathematical Sciences at Stony Brook
 for continued support and hospitality. We would also like to thank the
 audiences at the seminars in Stony Brook, Orsay and at the IHP
 where these
 results were presented.

\section{The Boundary of a Hyperbolic Component}
  \label{sec:boundary}

 \begin{defn}[Bifurcation Locus and Stable Components]
  The set 
   $\B := \{\kappa_0\in\C: \text{the family $\kappa\mapsto
          \Ek^n(\kappa)$ is not normal in $\kappa_0$}\}$
  is called the \emph{bifurcation locus} of exponential maps.
  A component $W$ of $\C\setminus\B$ is called
  \emph{hyperbolic} if all parameters in $W$ have an attracting
  periodic orbit. Otherwise,
  $W$ is called \emph{nonhyperbolic} or \emph{queer}.
 \end{defn}

  It is known that \emph{Misiurewicz parameters}, i.e. parameters for
  which the singular value is preperiodic, and 
  \emph{indifferent parameters}, that is,
  parameters with an indifferent periodic orbit, are dense in $\B$
  \cite{alexmisha,instability}.  

 Every attracting parameter lies in some hyperbolic
 component, and all
 parameters within the same hyperbolic component $W$ have the same
 period. As mentioned in the introduction, any hyperbolic component $W$ is
 simply connected and the multiplier map 
 $\mu:W\to\Ds$ is a universal covering
 map \cite[Theorem 12]{alexmisha}. Thus we can find a biholomorphic map $\Phi_W$ 
 from $W$ to the left halfplane $\H$ which satisfies
 $\exp\circ\Phi_W=\mu$. 
  We will often also consider the inverse
  $\Psi_W := \Phi^{-1}_W$ of this map. Note that
  $\Phi_W$ and $\Psi_W$ are uniquely defined
  only up to post- resp.\ precomposition by a translation of $2\pi i\Z$;
  this ambiguity will be removed in Proposition 
  \ref{prop:bifurcationstructure} by specifying
  a preferred choice of parametrization.

 Hyperbolic components in exponential parameter space do not
  have ``centers'': $\Ek$ has no critical points, so
  the attracting multiplier is never $0$. In some
  sense, their center is at $\infty$:

\begin{lem}[Small Multipliers] \label{lem:smallmultipliers}
 Let $W$ be a hyperbolic component of period $n$.
 If $(\kappa_i)$ is a sequence of parameters in a hyperbolic
 component $W$ such that $\mu(\kappa_j)\to 0$, then
 $\kappa_j\to\infty$. In particular, every hyperbolic component is unbounded.
\end{lem}
 \begin{proof}
   Let $\kappa\in W$ and let 
 $a_0\mapsto a_1\mapsto \dots\mapsto a_n=a_0$ be the attracting orbit
 of $\Ek$. 
 A simple calculation shows that
  \begin{equation}
    \mu(\kappa) = \prod_{j=1}^{n}
     \exp(a_j). \label{eqn:multiplieralongcycle}
  \end{equation}
 Since $|\mu(\kappa)|\leq 1$, there is some $j$ such that
 $\re a_j \leq 0$; without loss of generality let us suppose that
 $j=0$. This implies that
 $|a_1|\leq|\kappa|+1$. Denote $f(t) = \exp(t)+|\kappa|$; then
  $|a_k| \leq f^{k-1}(|\kappa|+1)\leq f^{n-1}(|\kappa|+1)$
 for $1\leq k\leq n$. 
 Thus the product (\ref{eqn:multiplieralongcycle}) is bounded from below
 in terms of $|\kappa|$, which concludes the proof.
\end{proof}

The hyperbolic components of
 periods one and two were described in
 \cite{bakerexp}; see Figure \ref{fig:hyperboliccomponents}. We will
 content ourselves here by noting that there is a unique component $W$
 of period one, which contains a left halfplane. Indeed, the map
  $\Psi_W$ is given by
  $\Psi_W(z)=z-\exp(z)$: for $\kappa=z-\exp(z)$, the
  point $z$ is a fixed point with multiplier $\mu=\exp(z)$.

We shall now turn our attention to points on the boundary of a hyperbolic component.

\begin{lem}[Indifferent Parameters]
  \label{lem:indifferentorbits}
  Let $W$ be a hyperbolic component of period $n$
   and $\kappa_0\in\partial W\cap\C$. Then
   $\kappa_0$ has an indifferent cycle of period dividing $n$.
   Furthermore, as
   $\kappa\to\kappa_0$ in $\cl{W}$, the nonrepelling cycles of $\Ek$ converge to the
   indifferent cycle of $E_{\kappa_0}$.
 \end{lem}
 \begin{proof}
   Let $\kappa_j\to\kappa_0$ in $\cl{W}$. By
  (\ref{eqn:multiplieralongcycle}), for every $j$ there exists some point
  $z_j$ on the
  nonrepelling orbit of $E_{\kappa_j}$ with $\re z_j \leq 0$.
  Thus the sequence $E_{\kappa_j}(z_j)=\exp(z_j)+\kappa_i$ is bounded and
  has some
  limit point $z\in\C$, which
  is
  a nonrepelling fixed point of $E_{\kappa_0}^n$.
  Since $E_{\kappa_0}$ has at most one nonrepelling orbit 
  \cite[Theorem 5]{alexmisha}, and since $z$ 
  cannot be attracting for $\kappa_0\in\partial W$,
   the claim follows.
 \end{proof}

Similarly to the definition of
external (and internal) rays for polynomials, the foliation of the
punctured disk by radial rays gives rise to a foliation of the
hyperbolic component by \emph{internal rays}. These rays are of a
natural interest when studying the boundaries of hyperbolic components.

\begin{defn}[Internal Rays] \label{defn:internalrays}
 Let $W$ be a hyperbolic component
 and let $h \in\R$. The curve
  \[ \IR:(-\infty,0)\to\C, t\mapsto \Psi_W( t + 2\pi i h)
  \]
 is called an \emph{internal ray at height $h$} (or at angle
 $\alpha$, where $\alpha$ is the fractional part of $h$).
 We say that an internal ray $\IR$
 \emph{lands} at a point $\kappa\in\Ch$ if $\kappa=\lim_{t\to 0} \IR(t)$.
\end{defn}
\begin{remark}
  By Lemma \ref{lem:smallmultipliers}, $\lim_{t\to -\infty} \IR(t)
   = \infty$. Note that the height of an internal ray is uniquely
   defined only if we fix a particular choice of $\Psi_W$; compare
   Proposition \ref{prop:bifurcationstructure} and the comment thereafter.
\end{remark}

\begin{lem}[Continuous Extension] \label{lem:localconnectivity}
 Let $W$ be a hyperbolic component. Then $\Psi_W$ extends to a continuous 
 surjection
 $\Psi_W:\cl{\H}\to\cl{W}$ with $\Psi_W(\infty)=\infty$. In particular,
 every internal ray of $W$ lands in $\Ch$.
\end{lem}
\begin{proof}
  It is sufficient to show that $\Psi_W$ has a continuous extension
 to $\cl{H}$; surjectivity then follows from the compactness of
 $\cl{H}$. 

 We must thus show that $\lim_{z\to z_0}\Psi_W(z)$
 exists for every $z_0\in i\R \cup \{\infty\}$. 
 Let $h\in\R$ and let $L$ denote the limit set of $\Psi_W(z)$ as
 $z\to ih$.
 By Lemma \ref{lem:indifferentorbits}, every parameter $\kappa_0\in L\cap\C$
 has an indifferent periodic point $a$ with multiplier $\mu = \exp(i
 h)$. We claim that the set of such parameters $\kappa$ is
 discrete. 

 Indeed, suppose there is a sequence $\kappa_j\in L$
 which converges nontrivially to some parameter
 $\kappa_0\in\C$. Let $a_0$ be an indifferent periodic point
 of $E_{\kappa_0}$; by
 Lemma \ref{lem:indifferentorbits}, we can pick a sequence $a_j\to a_0$
 such that each $a_j$ is on the indifferent periodic orbit
 of $E_{\kappa_j}$.
 The points $(\kappa_j,a_j)$ lie in the
 zero set of the function of two complex variables given by
  $f(\kappa,a) := \bigl(E_{\kappa}^n(a) - a , (E_{\kappa}^n)'(a) -
      \mu\bigr)$.
 Since $f$ is analytic, this implies that, for all $\kappa$
 in a neighborhood of $\kappa_0$, there is a solution of
 $f(\kappa,z)=0$. In other words, all parameters in a neighborhood of
 $\kappa_0$ are indifferent, which is absurd.

 Thus $L$ is contained in a totally disconnected set. Since $L$ is
 connected, this implies $\#L=1$ as required.

Finally, let us show that $\Psi_W(z)$ has no accumulation
 points in $\C$ as $z\to\infty$. Suppose then that $(z_j)$ is a
 sequence in $\H$ such that $\Psi_W(z_j)$
 converges to some
 point $\kappa_0\in\partial W\cap\C$.
 Then by Lemma \ref{lem:indifferentorbits},
 $\kappa_0$ has an indifferent orbit of multiplier, say, $e^{2\pi i
 h}$. It is easy to see that we can continue the multiplier of this orbit to
 an analytic function on a finite sheeted cover of a neighborhood
 $U$ of
 $\kappa_0$. In particular, there are finitely many connected subsets
 of $U$ in which the orbit becomes attracting, which means that (if
 $U$ was chosen small enough) the set $\Psi_W^{-1}(U)$ consists of
 finitely many bounded components. Since $z_j\in \Psi_W^{-1}(U)$ for
 large $j$, it follows that $z_j\not\to\infty$, as required. 
\end{proof}

\begin{cor}[Landing Points of Internal Rays] \label{cor:landingpoints}
 Let $W$ be a hyperbolic component.
 Then every point of $\partial W\cap\C$ is the landing point of a unique
 internal ray; in particular every component of $\partial W \cap \C$
 is a Jordan arc extending to $\infty$ in both directions.
\end{cor}
\begin{proof}
  The fact that every boundary point is the landing point of an
 internal ray follows immediately from Lemma
 \ref{lem:localconnectivity}. We need to show that no two
 internal rays can land at the same boundary point $\kappa_0\in\C$. Suppose that
  $\Psi_W(ih)=\Psi_W(ih')=\kappa_0$ for some $h<h'$. 

 Connect $ih$ and $ih'$ by a curve in $\H$; the image of
 this curve under $\Psi_W$ is then a simple closed curve $\gamma$
 in $\cl{W}$
 which intersects $\partial W$ only in $\kappa_0$. By the F.~and
 M.~Riesz theorem \cite[Theorem A.3]{jackdynamicsthird}, there exists some
 $h_1\in (h,h')$ with $\kappa_1:=\Psi_W(ih_1)\neq\kappa_0$.
 The indifferent parameter $\kappa_1$ lies in the
 bifurcation locus $\B$ and is separated from $\infty$ by the curve $\gamma$.
 Since Misiurewicz parameters are dense in $\B$, there exists
 some Misiurewicz parameter $\kappa_2$ which is also enclosed by
 $\gamma$. Indifferent parameters are also dense in $\B$, so we can
 find some indifferent parameter $\kappa$ with $|\kappa_2-\kappa| <
 \dist(\kappa_2,\partial W)$. The parameter $\kappa$ lies on the
 boundary of some hyperbolic component, which is thus
 separated from $\infty$ by $\gamma$. This
 contradicts Lemma \ref{lem:smallmultipliers}. 

 We have shown that $\Psi_W$ is injective on
 $\Psi_W^{-1}(\partial W \cap\C)$. It follows easily that
 every component $C$ of $\partial{W}\cap\C$ can be written as
 $C=\Psi_W(I)$, where $I$ is a component of $\Psi_W^{-1}(\partial
 W\cap\C)$. In other words, $I=i\cdot(h,h')$ with
 $\Psi_W(ih)=\Psi_W(ih')=\infty$, and thus $C$ is a Jordan
 arc tending to $\infty$ in both directions. 
\end{proof}

 We are now prepared to use the Squeezing Lemma to prove our main
 theorems.

\begin{proof}[Proof of Theorem {\ref{mainthm:boundary}}, using Theorem
 {\ref{thm:squeezing1}}.]
 We need to show that $\Psi_W:\cl{\H}\to\cl{W}$ is injective. By
 Corollary \ref{cor:landingpoints}, it only remains to show that  
 $\Psi_W(2\pi ih)\in\C$ for every
 $h\in\R$. 
 So suppose by contradiction that $h\in\R$ with $\Psi_W(2\pi
 ih)=\infty$, and
 define a curve $\gamma\colon[0,\infty)\to W$ by
 $\gamma(t):=\Gamma_{W,h}(-1/t)$. 

 Then $|\gamma(t)|\to\infty$ as 
 $t\to\infty$, while $\Phi_W(\gamma(t))=-1/t+2\pi ih\to 2\pi ih$. This 
 contradicts the Squeezing Lemma. 
\end{proof}

\begin{proof}[Proof of Corollary {\ref{cor:queer}} and
 Theorem {\ref{mainthm:bifurcationlocus}}, using Theorem
 \ref{thm:squeezing1}.]
  As already mentioned,
   indifferent and escaping parameters lie in the bifurcation locus, while
   attracting parameters lie in hyperbolic components. Hence it follows
   immediately from Theorem \ref{thm:squeezing1} that a non-hyperbolic
   stable component $U$ cannot contain a curve to infinity. 
   So $\infty$ is not an accessible point of $K = \Ch\setminus U$.
   It follows from the Plane Separation Theorem
   \cite[Theorem VI.3.1]{whyburnanalytictopology} that any cut-point of a
   non-degenerate full plane continuum $K$ is bi-accessible from $K$. (Compare
   \cite[Theorem 6.6]{mcmullenrenormalization}.) Hence it follows that
   $K\setminus\{\infty\} = \C\setminus U$ is connected. 
   This proves Corollary
   \ref{cor:queer}. 

  Finally, the bifurcation locus $\B$ could be disconnected only if
   some stable component $U$ disconnected $\C$. However, by Theorem
   \ref{mainthm:boundary} and Corollary \ref{cor:queer}, neither hyperbolic
   nor non-hyperbolic components can have this property, and Theorem
   \ref{mainthm:bifurcationlocus} is proved. 
\end{proof}

\section{Combinatorics of Exponential Maps} \label{sec:combinatorics}

 \subsection*{Dynamic rays}
 Let $\Ek$ be any exponential map. The \emph{set of escaping points}
 of $\Ek$ is defined to be
  \[ I := I(\Ek) := \{z\in\C: |\Ek^n(z)|\to\infty\}. \]
 It is known that the Julia set $J(\Ek)$ is the closure of $I(\Ek)$
 \cite{alexescaping,alexmisha}. Note that, because
 $|\Ek^n(z)|=|\exp(\re \Ek^{n-1}(z)) + \kappa|$, $z\in I(\Ek)$ if and only
 if $\re \Ek^n(z)\to+\infty$.

 A complete classification of the set of escaping points of an
 exponential map was given in \cite{expescaping}. To describe this
 result, let us introduce some combinatorial notation.
 An infinite sequence $\s = s_1 s_2 \dots$ of integers is called
 an \emph{external address}; we say that a point $z\in\C$ has external
 address $\s$ if 
  \[\im \Ek^n(z) \in \bigl( (2s_{n+1}-1)\pi , (2s_{n+1}+1)\pi\bigr) \]
 for all $n\geq 0$. An external address $\s$ is called
 \emph{exponentially bounded} if there exists some $x>0$ such that
  \[ 2\pi|s_n| < F^{n-1}(x) \]
 for all $n\geq 0$. (Recall that $F(t)=\exp(t)-1$.) 

\begin{thmdef}[Classification of Escaping Points \cite{expescaping}] \label{thm:classescaping}
 Let $\kappa\in\C$, and suppose that $\kappa\notin I(\Ek)$. 
 For every exponentially bounded address $\s$
 there exists $\ts\geq 0$ and a 
 curve $g_{\s}:= \gs^{\kappa}:(\ts,\infty)\to I(\Ek)$ or 
 $\gs := \gs^{\kappa}:[\ts,\infty)\to
 I(\Ek)$
 (called the \emph{dynamic ray} at address $\s$)
 with the following properties.
 \begin{enumerate}
   \item The trace of $\gs$ is a path connected component of $I(\Ek)$;
    \label{item:dynamicrayproperty1}
   \item for large $t$, $\gs(t)$ has external address $\s$;
   \item $\displaystyle{\Ek(\gs(t)) = g_{\sigma(\s)}(F(t))}$;
   \item $\gs(F^{n-1}(t))=F^{n-1}(t)+2\pi i s_n + O(e^{-F^{n-1}(t)})$ as $t$ or $n$
   tend to $\infty$.
    \label{item:dynamicrayproperty4}
 \end{enumerate}
 Conversely,
 every path connected component of $I(\Ek)$ is such a
 dynamic ray. 

 Now suppose the singular value does escape. Then there still exist
 dynamic rays $\gs$ with
 properties (\ref{item:dynamicrayproperty1}) to
  (\ref{item:dynamicrayproperty4})
  for all exponentially bounded addresses $\s$. 
 However, there are countably many $\s$ for which $\gs$ is
 not defined for all
 $t>\ts$ (resp.~$t\geq\ts$). More precisely, there exist $\s^0$ and
 $t^0\geq t_{\s^0}$ such that $\kappa=g_{\s^0}(t^0)$. For every
 $\s\in\sigma^{-n}(\s^0)$, the ray $\gs$ is not defined for $t\leq F^{-n}(t_0)$.
 Every path
 connected component of the escaping set
  is either a dynamic ray or is mapped into
 $g_{\s^0}$ by some forward iterate of $\Ek$. \qedd
\end{thmdef}

 As usual, we say that a ray $\gs$ \emph{lands} at a point $z\in\Ch$ 
   if $\lim_{t\to \ts} \gs(t)=z$.

 We need to be able to recognize a point which is on the
 dynamic ray $\gs$. Such a criterion is given by the following result,
 which is an extension of \cite[Theorem 4.4]{expescaping}. 

\begin{lem}[Fast Points are on Rays] \label{lem:fastpointsonray}
 Let $\kappa\in\C$, and let
  $x \geq \max(\re \kappa - 1,2\pi+6)$. 
  Suppose that $z_0\in I(\Ek)$ such that
   $\re \Ek^n(z_0) \geq F^n(x)$
  for all $n\geq 0$.
 Then $z_0 = \gs(t)$, where  $t\geq x$ and $\s$ is the external
  address of $z_0$. 
\end{lem}
\begin{remark} It is sufficient to require $x > Q(|\kappa|)$ with
  $Q(|\kappa|)=\log^+|\kappa| + O(1)$ (see \cite{topescapingnew}).
\end{remark}
\begin{sketch} By \cite[Theorem 4.3]{topescapingnew}, 
    $z_0=\gs(t)$ for some $t\geq \ts$. The fact
    $t\geq x$ follows by applying \cite[Theorem 4.4]{expescaping} to
    some forward iterate of $z_0$. 
\end{sketch}

\subsection*{Parameter rays}
 Given the importance of escaping points in providing structure in the
 dynamical plane, it is natural to ask about their analog in parameter
 space: \emph{escaping parameters}; i.e.~those 
 for which the singular value escapes.
 \begin{defn}[Parameter Rays] \label{defn:parrays}
  Let $\s$ be an exponentially bounded external address, let $t\geq
  \ts$. Then we define
    \[ \PR(t) := \{\kappa:\gs^{\kappa}(t)=\kappa\}. \]
  The set 
    $\PR := \bigcup_{t} \PR(t)$
  is called the \emph{parameter ray} at address $\s$.
 \end{defn} 
 In \cite{markus} (to appear in revised form as \cite{markusdierk})
 it was shown that, for $t>\ts$, the set $\PR(t)$
 consists of a single point, which depends continuously (and even
 differentiably) on $t$. Thus
 $\PR:(\ts,\infty)\to\C$ is a (differentiable)
 curve, justifying the term ``parameter ray''. 
 Using the Squeezing Lemma,
 this result is extended to escaping
 endpoints, and thus to a classification of all escaping parameters,
 in \cite{markuslassedierk}.

 For the purposes of this article, we will not require these
 results about parameter rays. However, we will sometimes use them heuristically
 in order to explain the strategy of a proof.

 \subsection*{Vertical order and intermediate external addresses}
  The existence of dynamic rays provides a structure to the Julia
  set and
  the dynamical plane. This structure can be used to describe the
  behavior of a 
  curve to $\infty$ under iteration, provided that the curve does not
  contain escaping parameters. Indeed, suppose that $\mathcal{C}$ is any disjoint family
  of curves $\gamma:[0,\infty)\to \C$ with $\re \gamma(t)\to+\infty$.
  Then $\mathcal{C}$ is
  equipped with a natural vertical order: of any two curves in
  $\mathcal{C}$, one is \emph{above} the other. More precisely, 
  define $\Hplane_R
  := \{z\in\C: \re z > R\}$ for $R>0$. If $\gamma\in\mathcal{C}$ and $R$
  is large enough, then  the set $\Hplane_R \setminus \gamma$ has exactly two
  unbounded components,
  one above and one below $\gamma$. Any other curve of $\mathcal{C}$
  must (eventually) tend to
  $\infty$ within one of these. 

  For the family of dynamic rays, it is easy to see that this vertical
  order
  coincides with the lexicographic order on
  external addresses.  We want to use this structure to assign combinatorics to 
  curves to $\infty$ in Fatou components. To do this, we 
  shall 
  add 
  \emph{intermediate external addresses} 
  to our repertoire. An intermediate external address is a finite 
  sequence 
  of the form 
   \[ \s = s_1 s_2 \dots s_{n-1} \infty, \] 
  where $n\geq 2$, $s_1,\dots,s_{n-2}\in\Z$ and $s_{n-1}\in 
  \Z+\frac{1}{2}$.  When we 
  wish to make the distinction, we will refer to an external address 
  in 
  the original sense as an ``infinite'' external address.
  For concise notation, we will always take the terms 
  ``exponentially bounded'', ``unbounded'' or ``bounded'' --- which 
  are not useful for intermediate addresses --- 
  to mean that the address is
  infinite and has the corresponding property.

  We 
  denote the space of all infinite and intermediate external addresses 
  by $\Sequ$. The point of this definition is that the space $\Sequ$ 
  is order-complete with respect to the lexicographic order; i.e., 
  every bounded subset has a supremum. We will also often consider 
  $\infty$ 
  as an intermediate external address of length $1$, and define 
  $\Sequb := \Sequ\cup\{\infty\}$. For more details, we refer the 
  reader to \cite[Section 2]{expcombinatorics}. 

 Now let $\gamma:[0,\infty)\to\C\setminus I(\Ek)$ be a curve to
  $\infty$. If $\lim_{t\to\infty} \gamma(t)=+\infty$, we define
  \begin{align*}
    \extaddr(\gamma) :=& \inf\{\s\in\Sequ: \text{$\s$ is exponentially
                                  bounded and $\gs$ is above $\gamma$}\} \\
                     =&  \sup\{\s\in\Sequ: \text{$\s$ is exponentially
                                  bounded and $\gs$ is below $\gamma$}\}.
  \end{align*}
  (Supremum and infimum exist by the
  order-completeness of $\Sequ$; they are equal since
  exponentially bounded addresses are dense in $\Sequ$.)
  If $\gamma$ does not tend to $+\infty$, then --- since $\gamma$ does
  not intersect dynamic rays ---
  $\re\gamma(t)$ is bounded from above. In this case, we set
  $\extaddr(\gamma)=\infty$. Note that, whenever $\extaddr(\gamma)\neq\infty$,
   $\extaddr(\Ek\circ\gamma)=\sigma(\extaddr(\gamma))$.

\subsection*{Attracting dynamics}
 Now suppose that $\Ek$ has an attracting or parabolic periodic
  cycle.   Then the singular value $\kappa$ is contained in
  some periodic Fatou component which we call the 
  \emph{characteristic Fatou component}.
  Let $U_0\mapsto U_1\overset{\approx}{\mapsto} 
              \dots \overset{\approx}{\mapsto} U_n=U_0$ 
  be the cycle of Fatou components, labeled such that
  $U_1$ is the characteristic component.
  Since $U_1$ contains a neighborhood of the singular value, $U_0$ contains an
  entire left half plane. In particular, $U_0$ contains a
  horizontal curve along which $\re(z)\to-\infty$, 
  which is
  unique up to homotopy. Its pullback $\gamma$ to $U_1$ under $\Ek^{n-1}$ is
  a
  curve to $+\infty$, and we define the \emph{intermediate external
  address of $\kappa$} to be
  $\extaddr(\kappa) := \extaddr(\gamma)$. Note that
  $\s:=\extaddr(\kappa)$ is an intermediate external address of length
  $n$ because $\extaddr(\Ek^{n-1}(\gamma))=\infty$.

It is easy to see that $\extaddr(\kappa)$ depends only on the
hyperbolic component $W$ which contains $\kappa$; this address will
therefore also be denoted by $\extaddr(W)$. (The same is true of other
combinatorial objects which we shall later associate to $\kappa$.) Its
significance lies in the following theorem, which is the main result
of \cite{expattracting}.
(The existence part of this theorem has also
appeared in \cite{dfj} and is a special case of Lemma \ref{lem:prescribedorbit}
below.)
\begin{thmdef}[Classification of Hyperbolic Components \cite{expattracting}]
 \label{thm:hyperboliccomponents}
 For every intermediate external address $\s$, there exists exactly one
 hyperbolic component $W$ with $\extaddr(W)=\s$. We denote this component
 by $\Hyp{\s}$. The vertical order of hyperbolic components coincides
 with the lexicographic order of their external addresses. \qedd
\end{thmdef}

To explain the last statement, let $\Phi_W:W\to \H$ with
$\exp\circ\Phi_W=\mu$, as before. The homotopy class of all curves
$\gamma:[0,\infty)\to W$ with $\Phi_W(\gamma(t))\to\infty$ is called
the \emph{preferred homotopy class of $W$}. As above, these preferred homotopy
classes have a vertical order, which is the order referred to
in the above theorem. 

\subsection*{A Precise Version of the Squeezing Lemma}

 Let $\gamma:[0,\infty)\to\C$ be a curve in parameter space. Suppose
  that $\lim_{t\to\infty} |\gamma(t)|=\infty$ and that $\gamma$ does not
  contain any indifferent parameters. Then, as above, we can
  associate to $\gamma$ an address
   \begin{align*}
     \extaddr(\gamma) := \inf\{\s\in&\Sequ: \s \text{ is an intermediate
      external address and} \\
      &\text{ every curve in the preferred homotopy
              class of $\Hyp{\s}$ is above
    $\gamma$}\}. 
   \end{align*}
  Note that, for every $k\in\Z+\onehalf$,
   the line $\{\im(\kappa)=2\pi k\}$ is contained in
   the union of
   the unique period $1$ component $\Hyp{\infty}$ (for $\re(\kappa)<1$),
   the period $2$ component $\Hyp{k\infty}$ (for $\re\kappa>1$)
   and their common
   parabolic boundary point $1+2\pi i k$. Therefore, if
   $\re(\gamma(t))\to+\infty$, then $\im(\gamma(t))$ is necessarily
   bounded, and $\extaddr(\gamma)$ starts with a finite entry. 
   In all other cases, $\re(\gamma(t))$ is bounded above and
   $\gamma$ is contained in $\Hyp{\infty}$; in this case,
   $\extaddr(\gamma)=\infty$.

  Using this definition, we can now state a more precise version of the
  Squeezing Lemma.

\begin{thm}[Squeezing Lemma] \label{thm:squeezinglemma}
 Let $\gamma:[0,\infty)\to\C$ be a curve in parameter space with
 $|\gamma(t)|\to\infty$. Suppose that $\gamma$
 contains no indifferent parameters. Then
 \begin{enumerate}
  \item[(a)] $\s := \extaddr(\gamma)$ is either intermediate or
    exponentially bounded.
  \item[(b)] If $\s$ is exponentially bounded, then  
              $\gamma(t)$ is escaping for all sufficiently
               large $t$. More precisely,
               $\gamma(t)\in\PR(\tau)$ for some
               $\tau\geq \re\gamma(t)-1$.
  \item[(c)] If $\s$ is intermediate, then $\gamma\subset \Hyp{\s}$ and
              $\lim_{t\to\infty}\Phi_{\Hyp{\s}}(\gamma(t))\to \infty$. (In
              other words, $\gamma$ lies in the
              preferred homotopy class of
              $\Hyp{\s}$.)
 \end{enumerate}
\end{thm}
Note that the proof of Theorem \ref{mainthm:boundary} in \cite{habil}
(outlined in \cite{cras}) consists of proving the Squeezing Lemma in
the case where $\gamma$ lies in some hyperbolic component. 
This is done by
distinguishing three cases, based on whether
$\s:=\extaddr(\gamma)$ is bounded, unbounded or intermediate. Our
proof for the exponentially bounded case (b) is a variant of the
bounded case in \cite{habil}: there exists a parameter ray at address
$\s$, and it is squeezed between nearby parameter rays (in
\cite{habil}) resp.\ by hyperbolic components (here).
In the intermediate case, both proofs use the bifurcation structure of
hyperbolic components, showing that there are curves which separate
$\gamma$ from the preferred homotopy class of 
$\Hyp{\s}$. However, in our case --- where we do not know that $\gamma$ 
starts out in some hyperbolic component --- we need an additional 
tool (Lemma \ref{lem:parameterraybound} below) 
to ensure that the real parts on these curves tend to 
$+\infty$. 

It remains to exclude the unbounded resp.\ exponentially unbounded 
cases. Both proofs use the concept of 
\emph{internal addresses} to show that $\s$ is contained in 
the nested intersection of 
wakes of infinitely many hyperbolic components such that this 
intersection contains only $\s$, 
and hence no 
hyperbolic components. This
suffices to exclude the case of unbounded
addresses for curves in hyperbolic components.
For our proof, we again use Lemma
\ref{lem:parameterraybound} to show that, in the exponentially
unbounded case, these
wakes have real parts tending to $+\infty$, and thus their
intersection is empty.

\section{Squeezing around Parameter Rays}
\label{sec:squeezingb}

In this section, we shall prove part (b) of the Squeezing Lemma, which
we restate here for convenience.

\begin{figure}
   \center
  \subfigure[Idea of the proof of Theorem \ref{thm:squeezingb}]{
   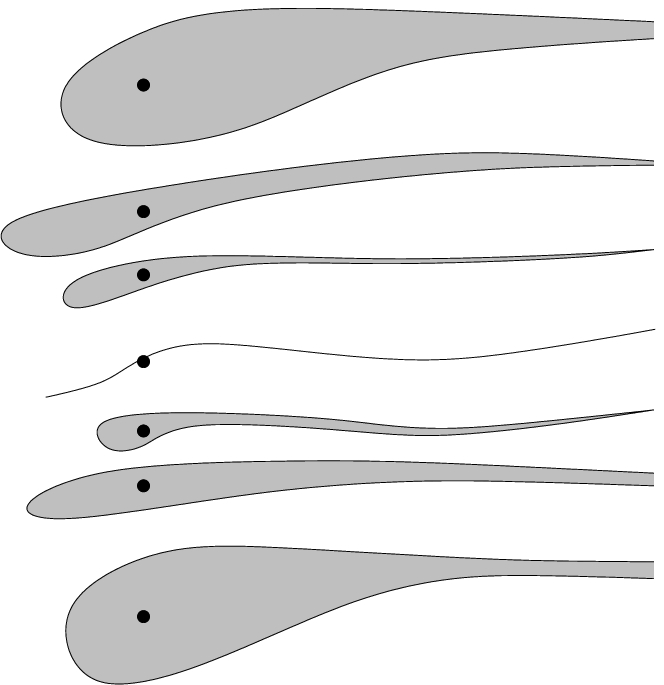
      \label{fig:squeezingb_marked}} \hfill
  \subfigure[Illustration of Lemma \ref{lem:prescribedorbit}]{
    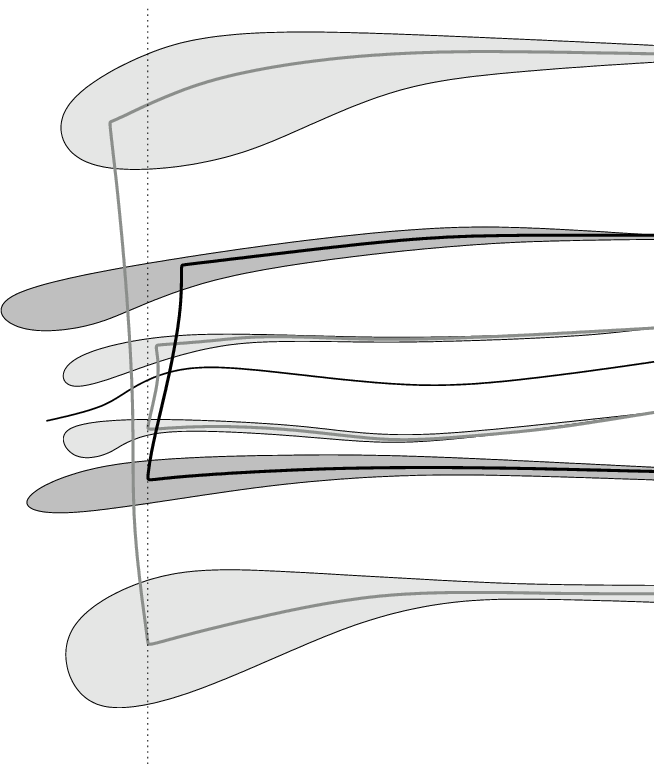
      \label{fig:squeezingb_strips}}
  \caption{Squeezing along parameter rays}
\end{figure}

\begin{thm}[Curves at Exponentially Bounded Addresses] \label{thm:squeezingb}
 Let $\gamma:[0,\infty)\to\C$ be a curve in parameter space which does
 not contain any indifferent parameters. Suppose that
 $\gamma(t)\to\infty$ and that $\s := \extaddr(\gamma)$ is exponentially
 bounded.
 Then, if $t$ is large enough, $\gamma(t)$ is escaping; in fact
 $\gamma(t)\in\PR(\tau)$ with $\tau\to\infty$ as $t\to\infty$.
\end{thm}

This fact was originally proved --- for bounded external addresses ---
by approximating the parameter ray at a given address 
above and
below by other parameter rays. This
is the origin of the name ``Squeezing Lemma''. Our proof will instead
squeeze from above and below using hyperbolic components. It is a generalization
of the proof of existence of
hyperbolic components \cite[Theorem 3.5]{expattracting}; the essence of the argument
goes back as far as \cite[Section 7]{bakerexp} (which, however, did
not use a combinatorial description to distinguish these components).

Before we come to the technical details of the proof, let us outline
 the general idea. Let us suppose for a moment that we know that the
 ``parameter ray'' $\PR$
 is indeed a ray with
 $\extaddr(\PR)=\s$. When we pick such a parameter $\kappa_0$
 sufficiently far to the right, the singular orbit will escape to
 infinity exponentially fast and with external address $\s$. Now, for
 every sufficiently large $n$, we
 should be able to perturb $\kappa_0$ slightly to parameters
 $\kappa_n^{\pm}$
 so that the $n$-th image
 of the singular value moves to the lines $\{\im z = 2\pi
 (s_{n+1}\pm\onehalf)\}$, and thus is mapped (far) to the left in the next
 step. The parameters $\kappa_n^{\pm}$ should then have attracting orbits
 and intermediate external addresses 
 $s_1 \dots s_{n} (s_{n+1}\pm\onehalf)\infty$, and
 the size of the perturbation of $\kappa_0$
 will tend to zero as $n\to\infty$. The
 perturbed parameters can be connected to $\infty$ within their
 hyperbolic components, and thus the
 curve $\gamma$ has to pass in between every pair
 $(\kappa_n^-,\kappa_n^+)$, and therefore contains $\kappa_0$. (Compare
 Figure \ref{fig:squeezingb_marked}.)

In order to carry out this plan, we need to be able to identify
 parameters in a given hyperbolic component. 
 We shall call a parameter $\kappa\in\C$ \emph{exponentially $\kappa$-bounded for
 $n$ steps} if $\re\kappa \geq 4$, 
 $|\im(\Ek^{k-1}(\kappa))|<F^{k-1}(\re\kappa-2)$ for 
 $1\leq k \leq n$ and $\re(\Ek^{k-1}(\kappa))>0$ for
 $1\leq k\leq n-1$. We shall also say that $\kappa$ has \emph{initial
 external address $s_1 \dots s_n$} if 
 $\im(\Ek^{k-1}(\kappa))\in \bigl((2s_k-1)\pi , (2s_k+1)\pi\bigr)$ for
 $1\leq k \leq n$.

\begin{lem}[Exponential Growth] \label{lem:expogrowth}
 Let $\kappa\in\C$ such that $\re\kappa\geq 2$. 
  If $z\in\C$  with $\re z \leq \re\kappa + 1$, then
  \begin{equation}
    |\re \Ek^{k}(z)| \leq F^{k}(\re(\kappa)+2)
       \quad\text{and}\quad
    |\Ek^{k}(z)| \leq F^{k}(|\kappa|+2)  \label{eqn:boundonorbit}
  \end{equation} 
 for all $k\geq 1$. 
 Furthermore, if $\kappa$ is exponentially
 $\kappa$-bounded for $n\geq 1$ steps, then
 \begin{equation}
    |\re(\Ek^{k}(\kappa))| \geq F^{k}(\re(\kappa)-1)
     \quad \text{for $0\leq k\leq n-1$}.
                          \label{eqn:largerealpart} 
 \end{equation} 
\end{lem}
\begin{proof}
  First note that, for $t,x\geq 2$,
  \begin{align}
    F(t) &\geq 2t + 2 \label{eqn:exponentialinequality1}\\
    F(x+1) &\geq 2F(x)+1 \geq F(x)+2F(x-1)+1.  \label{eqn:exponentialinequality2}
  \end{align}  
 To prove the first claim, let us show that the slightly stronger
  inequality 
  \begin{equation}
   |\re\Ek^k(z)| \leq F^k(\re\kappa+2)-1 \label{eqn:inductiveclaim}
  \end{equation}
  holds for all $k\geq 0$. Indeed it holds
  for $k=0$ by assumption, and, inductively,
   \begin{align}
   \notag
       |\re \Ek^{k+1}(z)| + 1 &\leq 
       \exp(\re\Ek^k(z))+\re\kappa+1 = F(\re\Ek^{k}(z)) + 2 + \re \kappa 
               \\ \notag 
  &\leq
       F(F^k(\re\kappa+2)-1) + 2 + \re\kappa \\ \notag
  &\underset{(\ref{eqn:exponentialinequality1})}{<}
      2F(F^k(\re\kappa+2)-1) 
      \underset{(\ref{eqn:exponentialinequality2})}{<}
            F^{k+1}(\re\kappa+2).
   \end{align}
  The second inequality regarding absolute values
   follows from (\ref{eqn:inductiveclaim}) by similarly noting that
  \[ |\Ek^{k+1}(z)|\leq F(\re\Ek^k(z))+1+|\kappa| \leq
       F(F^k(\re\kappa+2)-1)+1+|\kappa| < F^{k+1}(|\kappa|+2). \]

 Now suppose $\kappa$ is exponentially $\kappa$-bounded for $n$ steps. 
  We set $x:=\re\kappa -1$, $t_0:=x-1$ and $z_k := \Ek^{k-1}(\kappa)$. 
  In order to
  prove (\ref{eqn:largerealpart}), we 
  prove inductively that
  $|\re( z_{k} )| \geq F^{k-1}(x) + 1$ for $1\leq k\leq n$. 
 Indeed, suppose that $k\geq 0$ such that this claim is true for
  $k\geq 1$. 
  Then $|\re(z_{k+1})|\geq |z_{k+1}|- F^k(t_0)$ by hypothesis, 
   and we can estimate
   \begin{align*}
       |\re(z_{k+1})| & \geq |z_{k+1}|-F^k(t_0)
    \geq
        \exp(\re(z_k)) - |\kappa| - F^k(t_0)  \\
    &\geq F(\re(z_k))-\re(\kappa) - |\im(\kappa)| - F^k(t_0) \\
    &\geq F(F^{k-1}(x)+1)- (2t_0+2) - F^k(t_0) \\
    &\underset{(\ref{eqn:exponentialinequality1})}{\geq}
                   F(F^{k-1}(x)+1) - 2F^k(t_0) \\
    &\geq F(F^{k-1}(x)+1) - 2F(F^{k-1}(x)-1) 
     \underset{(\ref{eqn:exponentialinequality2})}{\geq} F^k(x) + 1. 
   \end{align*}
 (In the fourth and fifth step, we used the facts that
  $|\im\kappa|\leq \re\kappa-2$ and $t_0\geq \re\kappa\geq 2$,
  respectively.)
   \end{proof}

The following lemma generalizes
 \cite[Lemma 3.4]{expattracting}.
\begin{lem}[Identifying Hyperbolic Components] \label{lem:paramderiv}
 Suppose that $\kappa_0\in\C$ is exponentially
 $\kappa_0$-bounded for $n\geq 2$ steps, and consider the map
 $f_n:\kappa\mapsto \Ek^{n-1}(\kappa)$. Then 
 \begin{equation}
    |f_n'(\kappa_0)| \geq
                    F^{n-1}(\re(\kappa_0)-1). \label{eqn:paramderiv} 
 \end{equation}

 If $\re(E_{\kappa_0}^{n-1}(\kappa_0))<0$, then 
  $E_{\kappa_0}$ has an attracting
  periodic orbit of exact period $n$. The multiplier of this orbit
  tends to $0$ as either $\re(\kappa_0)$ or $n$ tend to
  $\infty$. If $s_1 \dots s_{n-2}$ is the initial
  external address of $\kappa_0$ and $s_{n-1}\in\Z+\onehalf$ with 
  $\im E_{\kappa_0}^{n-1}(\kappa_0) \in \bigl((2s_{n-1}-1)\pi,(2s_{n-1}+1)\pi\bigr)$,
  then $\extaddr{\kappa_0}=s_1 \dots
  s_{n-1}\infty$.
\end{lem}
\begin{proof}
  Let us again abbreviate $z_k := E_{\kappa_0}^{k-1}(\kappa_0)$.
 We shall show  (\ref{eqn:paramderiv}) by induction in $n$. So suppose
  that either $n=2$, or that $n>2$ and (\ref{eqn:paramderiv}) is true
  for $n-1$. By the chain rule,
  \[ f_n'(\kappa_0) = 
      \exp(z_{n-1})f_{n-1}'(\kappa_0) + 1.
   \]
 Furthermore, $|f_{n-1}'(\kappa_0)|\geq 1$ either trivially (if $n=2$) or by
 the induction hypothesis (if $n>2$). 
 Thus (\ref{eqn:largerealpart}) implies that
  \[ |f_n'(\kappa_0)| \geq |\exp(z_{n-1})|-1 =
           F(\re(z_{n-1})) \geq 
           F^{n-1}(\re(\kappa_0)-1). \]

 Now suppose that $\re z_n<0$. By Lemma \ref{lem:expogrowth}, 
  \[ \re z_{n} < -F^{n-1}(\re\kappa_0 - 1 ). \]  
  The fact that $E_{\kappa_0}$ is attracting (with exact period $n$) was
  proved in \cite[Lemma 3.4]{expattracting}
  under the stronger assumption that $\im z_k$ is bounded. The proof
  remains essentially the same,
  so we shall only sketch it without working out
  the precise estimates. The image of the left
  half plane $\Hplane:=\{z\in\C: \re z < \re z_n +1 \}$ under $E_{\kappa_0}$ is a
  punctured disk
  $D$ around $\kappa_0$ with
  radius at most 
    $\exp(1-F^{n-1}(\re\kappa_0-1))$.
  Using (\ref{eqn:boundonorbit}), it is not
  difficult to show that, for any point $z\in D$,
    \[ \Bigl|\bigl(E_{\kappa_0}^{n-1}(z)\bigr)'\Bigr|=
      \prod_{k=1}^{n-1} \exp(\re(E_{\kappa_0}^{k-1}(z))) 
      \leq
        \exp(F^{n-1}(\re\kappa_0-2)). \]
  Let $U := E_{\kappa_0}^{n-1}(D)=E_{\kappa_0}^{n}(\Hplane)$; it follows that
   \[ \diam U \leq \exp(1 + F^{n-1}(\re\kappa_0 -2) -
       F^{n-1}(\re(\kappa_0)-1))\leq 1 \]
  Thus $U\Subset \Hplane$, and $U$ contains an attracting periodic point $a$
   of exact period $n+1$.  As $n$ or
  $\re\kappa_0$ gets large, the diameter of $U$ tends to $0$. Since
  $U$ contains the image of $\D_1(z_n)$, the multiplier
  $\mu=\bigl(E_{\kappa_0}^{n-1}\bigr)'(a)$ 
   must also tend to $0$ by Koebe's theorem. 

 To see that $\kappa_0$ has the correct external address, connect
 $z_{n}$ to $-\infty$ by a curve $\gamma_{n}$ at constant imaginary parts,
 and consider the pullbacks $\gamma_k$ of this curve along the orbit
 $(z_k)$. We claim that, for $1<k\leq n-1$, all points on $\gamma_k$ have
 real parts larger than $\re \kappa_0$. (This shows that these
 curves do not cross the partition boundaries, and therefore
 completes the proof.) Indeed, 
 suppose by contradiction that $z\in \gamma_k$ with $\re z \leq
 \re\kappa_0$. Then,
  by (\ref{eqn:boundonorbit}),
 $|\re E_{\kappa_0}^{n-k}(z)| \leq F^{n-k}(\re\kappa_0+2)$, which contradicts
 the fact that
 \[ |\re E_{\kappa_0}^{n-k}(z)|\geq |\re z_{n}| 
   \geq F^{n-1}(\re\kappa_0-1) 
   \underset{(\ref{eqn:exponentialinequality1})}{\geq}
   F^{n-2}(2\re\kappa_0) > 
   F^{n-2}(\re\kappa_0+2). \qedhere\] 
\end{proof}

 For any infinite external address $\s$ and any $n\geq 2$, let us define
  \[ \s^{n}_{\pm} :=  s_1 s_2 \dots s_{n-2} (s_{n-1}\pm\onehalf) \infty,  \] 
 and consider, for any $R>0$, the half strip
  \[ K_{n,R}:=\{z\in\C:\re(z)\geq R, \im(z)\in [(2 s_{n-1}-1)\pi,
       (2s_{n-1}+1)\pi]\}. \]
 The statement of the next lemma is somewhat technical, but the
  general idea is very simple.
  Our goal is to construct, for every $n$, an inverse branch
  $\PBM$ of the map $f_n:\kappa\mapsto \Ek^{n-1}(\kappa)$, defined on a suitable
  $K_{n,R_n}$, such that the upper boundary of $K_{n,R_n}$ is mapped into
  $\Hyp{\s^n_+}$ and the lower into $\Hyp{\s^n_-}$. A curve $\gamma$ as in
  the statement of Theorem \ref{thm:squeezingb} will
  then have to (eventually) lie in the image of every $\PBM$, and
  as $n$ grows,
  these will be thinner and thinner by the derivative estimate
  (\ref{eqn:paramderiv}) of Lemma \ref{lem:paramderiv}. 
  (Compare Figure \ref{fig:squeezingb_strips}.)

\begin{lem}[Parameters with Prescribed Singular Orbit]
  \label{lem:prescribedorbit}
 Let $\s$ be any exponentially bounded address, say
  $2\pi(|s_k|+1)<F^{k-1}(x)$ with $x\geq 3$. Set $R:=x+2$ and
  let $n\geq 2$ be arbitrary. 
 Then
 there is $R_n$ with $|F^{-(n-2)}(R_n)-R|\leq 1$
 and a conformal map $\PBM$ from
 $K_n:=K_{n,R_n}$
 into parameter space such that
 \begin{enumerate}
  \item[1.] For any $z\in K_n$, $\kappa:=\PBM(z)$ is
        exponentially $\kappa$-bounded
        for $n-1$ steps and has initial external address 
        $s_1 \dots s_{n-2}$. 
  \item[2.] $f_{n-1}(\PBM(z))=z$ for all $z$. 
  \item[3.] $\re\Bigl(\PBM\bigl(R_n+(2s_{n-1}-1)\pi i\bigr)\Bigr)=R$.
 \end{enumerate}
\end{lem}
\begin{remark} This proves the existence of hyperbolic
 components with an arbitrary intermediate external address; i.e.~one
 half of Theorem \ref{thm:hyperboliccomponents}.
\end{remark}
\begin{proof}
 For all $\kappa$,
  set $T_{\kappa}:=\max\{x+1,\log(2(|\kappa|+2))\}$. A simple induction
  (quite similar to that
  of \cite[Lemma 3.3]{expescaping} or \cite[Lemma 4.1]{topescapingnew}) 
   shows that there exists a branch 
   $\phi_{\kappa}$ of $\Ek^{-(n-2)}$ on $K_{n,F^{n-2}(T_{\kappa})}$
  such that, for all $z\in K_{n,F^{n-2}(T_{\kappa})}$,
   \begin{itemize}
    \item
     $\phi_{\kappa}(z)$ has initial external address $s_1\dots s_{n-2}$
        and
    \item $|\re(\Ek^{k-1}(\phi_{\kappa}(z)))-F^{k-n+1}(\re z)|<1$ for
      $1\leq k \leq n-1$.
   \end{itemize}
  Observe that $\phi$ depends holomorphically on $z$ and $\kappa$.

We claim that there exists $\kappa_0$ with $\re\kappa_0=R$
  such that
  $\kappa_0=\phi_{\kappa_0}(r_0+(2s_{n-1}-1)\pi i)$ for some $r_0\geq
  T_{\kappa_0}$. 
  Indeed, note that, if
  $\re \kappa = R$ and 
  $\im\kappa\in \bigl( (2s_1-1)\pi , (2s_1+1)\pi \bigr)$, then
  $T_{\kappa} = x = R-2$. Thus the complement of the set 
   \[ A_{\kappa}:=\{z\in\C:\re z \leq R-1\}\cup
                  \{\phi(r+(2s_{n-1}-1)\pi i):r\geq T_{\kappa}\} \]
  has two unbounded components, one above and one below $A_{\kappa}$. 
  For $\kappa=R+(2s_1-1)\pi i $, the singular
  value lies in the lower of these components, and for
  $\kappa=R+(2s_1+1)\pi i$ it lies in the upper. Thus, for an
  intermediate choice $\kappa_0$, the singular value must be contained
  in $A_{\kappa_0}$.

 So let $\kappa_0$ and $r_0$ be as above, and 
  set $R_n:= r_0$, $K_n:=K_{n,R_n}$ 
  and $\PBM(R_n+2(s_{n-1}-1)\pi i):=\kappa_0$.
  We claim that we can extend $\PBM$ to an analytic
  function $\PBM:K_n\to\C$ with $\phi_{\PBM(z)}(z)=\PBM(z)$. 
  First note that, whenever $\phi_{\kappa}(z)=\kappa$, the parameter
  $\kappa$ is exponentially $\kappa$-bounded for $n-1$ steps and has initial
  external address $s_1 \dots s_{n-2}$. Thus, by
  (\ref{eqn:paramderiv}) and the implicit function theorem, we can
  locally extend any such solution to an analytic function of $z$.

 Since $K_n$ is simply connected, by the monodromy theorem it
 suffices to show that we can continue $\PBM$ analytically along every
  curve $\gamma:[0,1]\to K_n$ with $\gamma(0)=r_0+2(s_{n-1}-1)\pi i$. Let
  $I$ be the maximum interval such that the solution can be continued
  along $\gamma|_I$. By the above remark, $I$ is open as a subset of
  $[0,1]$. It thus remains to show that $I$ is closed. So let $t_0$ be
  a limit point of $I$, and let 
  $\kappa$ be a limit point of $\PBM(\gamma(t))$ as $t\to 
  t_0$. Then, by continuity of $\phi$, 
  $\phi_{\kappa}(\gamma(t_0))=\kappa$. By 
  (\ref{eqn:paramderiv}), the set of such $\kappa$ is discrete, and 
  thus $\PBM(\gamma(t))\to\kappa$ as $t\to t_0$.  \end{proof}

\begin{proof}[Proof of Theorem \ref{thm:squeezingb}.] 
 Set $\s:=\extaddr(\gamma)$, 
 and choose $x\geq 3$ with $2\pi (|s_k|+1)\leq F^{k-1}(x)$. For 
  every 
  $n$, let $\PBM$, $R_n$ and $K_n$ be as in Lemma 
  \ref{lem:prescribedorbit}.

  By Lemma \ref{lem:paramderiv}, the curves 
    $\gamma^{n}_{\pm}:t\mapsto \PBM(t+2(s_{n-1}\pm 1)\pi i)$ lie in
    the preferred homotopy
    classes of $\Hyp{\s^n_{\pm}}$. Note also that, by the estimate
  (\ref{eqn:paramderiv}), 
  \[ \diam \PBM\Bigl(\bigl\{R_n + bi: b \in
           [(2s_{n-1}-1)\pi,(2s_{n-1}+1)\pi]\bigr\}  \Bigr) 
           \leq \frac{2\pi}{F^{n-2}(x+1)}. \]
 Recall that $\gamma$ tends to $\infty$ between the components $\Hyp{\s^n_-}$
  and $\Hyp{\s^n_+}$, and thus between $\gamma^{n}_-$ and $\gamma^n_+$.
  It follows that, for large $t$,
    \[ \gamma(t) \in \bigcap_n \PBM(K_n). \]
  By Lemmas \ref{lem:expogrowth} and \ref{lem:prescribedorbit}, this
  means that $\gamma(t)$ is escaping, and the singular orbit of
    $E_{\gamma(t)}$ has external address
  $\s$. By (\ref{eqn:largerealpart}) and Lemma
  \ref{lem:fastpointsonray}, $\gamma(t)=\gs^{\gamma(t)}(\tau)$ with 
  $\tau\geq \re\gamma(t)-1$. 
\end{proof}


\section{Bifurcations of Hyperbolic Components} \label{sec:wakes}
 In this section, we shall collect the combinatorial notions and
 results from \cite{expattracting}, \cite{expper}
 and \cite{expcombinatorics} which we
 require to complete the proof of the Squeezing Lemma.

\subsection*{Itineraries and Kneading Sequences}

 An important aspect of attracting (as well as escaping and singularly
 preperiodic) exponential dynamics is the possibility to associate so-called
 \emph{itineraries} to orbits of $\Ek$. The idea is that, whenever
 there is a curve (such as a dynamic ray or a curve in a Fatou
 component) starting at $\infty$ and 
 ending at the singular value, the preimages of this curve
 will form a partition of the plane. Here we will only
 define the combinatorial analog of this notion and one of its
 important 
 applications; compare \cite{expper}
 for more details and motivation.

 \begin{defn}[(Combinatorial) Itinerary] \label{defn:combinatorialitinerary}
  Let $\adds\in \Sequ$ and $\r\in\Sequb$. 
  Then
   the \emph{itinerary of $\r$ with respect to $\adds$} is
   $\itin_{\adds}(\r)=\u_1 \u_2 \dots$, where
   \[
       \begin{cases}
        \u_k = \j      & \text{if \hspace{1mm}
                                 $\j\adds < \sigma^{k-1}(\r) < (\j+1)\adds$} \\
        \u_k = \itj &\text{if \hspace{1mm}
                                 $\sigma^{k-1}(\r)= \j\adds$}  \\
        \u_k = {\tt *}       & \text{if \hspace{1mm}
                                 $\sigma^{k-1}(\r)=\infty$}.
      \end{cases}\]
 \end{defn}
 \begin{remark} Thus $\itin_{\s}(\r)$ has the same length as $\r$; in
  particular, itineraries of intermediate external addresses are finite.
 \end{remark}

 \begin{prop}[Rays Landing Together {\cite[Proposition 4.5]{expper}}]
  Let $\kappa$ be an attracting or parabolic parameter and
  $\s:=\extaddr(\kappa)$. Let $\r$ and $\rt$ be periodic external
  addresses. Then the dynamic rays
  $g_{\r}$ and $g_{\rt}$ land at a common point if and only if
  $\itin_{\s}(\r)=\itin_{\s}(\rt)$.  \qedd
 \end{prop}
 \begin{remark}
  This result remains true without the assumption of
   periodicity; see \cite{tying,topescapingnew}.
 \end{remark}

 As in the case of external addresses, the itinerary which describes the
  behavior of the singular value is of particular importance:

 \begin{defn}[Kneading Sequence]
  Let $\s\in\Sequ$. Then the \emph{kneading sequence} of $\s$ is defined
   to be $\K(\s) := \itin_{\s}(\s)$; we also define 
   $\K(\infty) := {\tt *}$.
  If $\s$ is an intermediate external address and $W=\Hyp{\s}$ is the
  corresponding hyperbolic component, we also set $\K(W) := \K(\s)$.
 \end{defn}

\subsection*{Characteristic Rays and Wakes} 
 If a quadratic polynomial has an attracting orbit, then the combinatorics
  (and hence the topology) of the Julia set can be described exactly
  in terms 
  of the periodic
  rays which land together on the boundaries of
  periodic Fatou components. More
  precisely, if $p$ is a hyperbolic quadratic polynomial, then
  there is a unique pair of periodic rays, called the
  \emph{characteristic rays}, landing together on the
  boundary of the Fatou component containing the critical value and
  separating it from all other points on the attracting orbit. The
  corresponding parameter rays land together at the hyperbolic
  component containing $p$ and bound the \emph{wake} of this
  component.

\begin{figure}
 \subfigure[]%
  {\resizebox{.48\textwidth}{!}{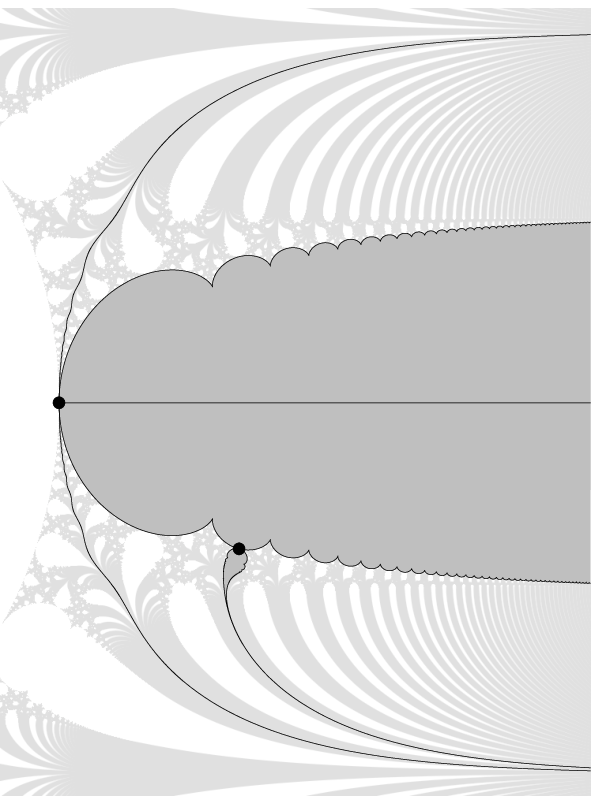\label{fig:bifurcationstructure}}}%
 \hfill%
 \subfigure[]%
  {\resizebox{.48\textwidth}{!}{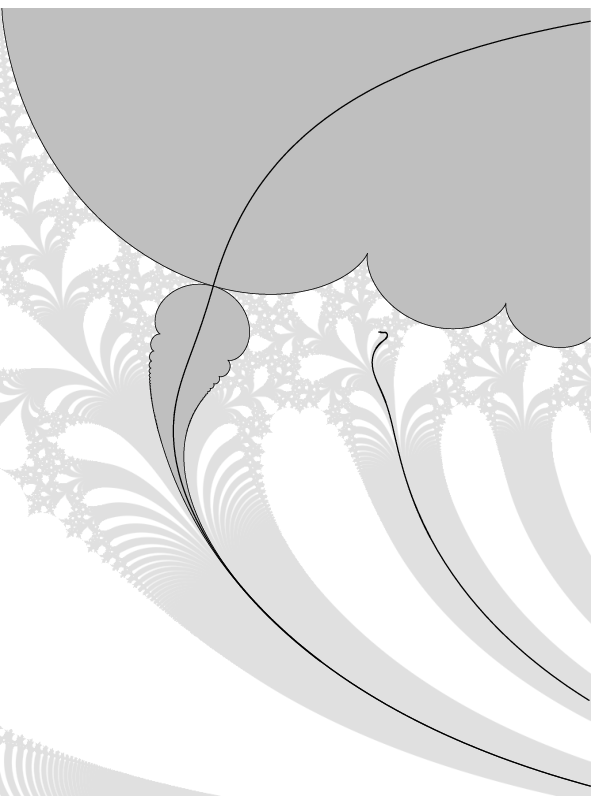\label{fig:cuttinginwake}}}%
 \caption{Illustration of Propositions \ref{prop:bifurcationstructure}
  and \ref{prop:cuttinginwake} for the component
  $W=\Hyp{\onehalf\infty}$. In (a), the two parameter rays
  bounding the wake of $W$ are shown; also drawn in is the child
  component $\child$ with $h=3/2$. (b) exemplifies Proposition
  \ref{prop:cuttinginwake}; here the child component
  is $V:=\Bif(W,\onehalf)=\Hyp{01\onehalf\infty}$. The curve
  running through these two components surrounds $\gamma$,
  which in this case is the parameter ray at address
  $\per{02030}$.}
\end{figure}

 The analogous dynamical statements were shown for exponential maps in
  \cite{expattracting}; they will be used below to define a
  combinatorial notion of wakes. In fact, it is also true (see Section
  \ref{sec:furtherresults}) in the
  exponential family that, on the boundary of every hyperbolic
  component, there are two parameter rays landing together which
  provide a natural definition of the \emph{wake} of this component
  as a subset of exponential parameter space 
  (Figure \ref{fig:bifurcationstructure}). While we
  can and will not use this result in our proofs, it often provides a
  useful motivation for the combinatorial ideas.

 \begin{prop}[Characteristic Rays \protect{\cite[Theorem 6.2]{expattracting}}]
   \label{prop:characteristicrays}
  Let $\s$ be an intermediate external address of length $n\geq 2$,
  and let $W=\Hyp{\s}$ be the hyperbolic component of address $\s$. Then there exists a
  (unique) pair of periodic external addresses $\s^-$ and $\s^+$ of period
  $n$ with the following properties.
  \begin{itemize}
   \item $\s^-<\s<\s^+$;
   \item $\itin_{\s}(\s^-)=\itin_{\s}(\s^+)$, and the first $n-1$ entries
           of this common itinerary agree with $\K(\s)$;
   \item $\sigma^j(\s^-),\sigma^j(\s^+)\notin (\s^-,\s^+)$ for all $j$
         (where $(\s^-,\s^+)=\{\r\in\Sequ: \s^-<\r<\s^+\}$).
  \end{itemize}
  These addresses are called the \emph{characteristic external
  addresses} of $\s$ (or $W$). Their (common) itinerary is called the \emph{forbidden
  kneading sequence of $W$} and denoted by $\KS(W)$. The interval
  $\W(W) := (\s^-,\s^+)\subset\Sequ$ is called the
  \emph{(combinatorial) wake} of $W$. \qedd
 \end{prop}

 \begin{prop}[Characteristic Rays and Wakes {\cite[Lemma 3.9]{expcombinatorics}}] 
  \label{prop:wake}
  Let $W$ be any hyperbolic component, and let $\s^-$ and $\s^+$ be
  its characteristic external addresses.
  Then, for any $\r\in\Sequ$,
  $\itin_{\r}(\s^+)=\itin_{\r}(\s^-)$ if and only if $\r\in\W(W)$. \qedd
 \end{prop}

\subsection*{Bifurcation Structure} Finally, we will need a good
  description of the structure of hyperbolic components; particularly
  how these components are connected to each other. Again, suitable
  analogs of many of
  the following statements are well-known for Mandel- and
  Multibrot sets; see e.g.~\cite{jackrays,fibers}.

 Recall that, for any hyperbolic component $W$, there exists a
  conformal isomorphism $\Psi_W:\H\to W$ with $\mu\circ\Psi_W=\exp$, 
  which extends continuously to $\cl{\H}$ by Lemma
  \ref{lem:localconnectivity}. Recall also that this defines $\Psi_W$ 
  uniquely up to precomposition by an 
  additive translation of $2\pi i\Z$.
  We will now
  describe the structure of the hyperbolic components
  which touch $W$ in its parabolic boundary points, using a preferred
  choice for $\Psi_W$. 

 The statement we will utilize (Proposition
  \ref{prop:bifurcationstructure} below) is somewhat abstract, so let
  us outline the most important ideas. Suppose that
  $h=\frac{p}{q}\in\Q\setminus\Z$ with
  $\kappa_0:=\Psi_W(2\pi i h)\in\C$. Then $\kappa_0$ is a parabolic
  parameter, and by arguments similar to those
  used in the Mandelbrot set \cite{jackrays}, there exists a unique
  component $\child$ of period $qn$ which touches $W$ at
  $\kappa_0$ (here $n$ is the period of $W$). Let us call the address
  of this component $\extaddr(W,h):=\extaddr(\child)$; this
  address lies in the wake of $W$, either
  above or below $\s$.  As $h$ increases, there will be some point
  $h_0$ at which $\extaddr(W,h)$ jumps from above $\s$ to
  below $\s$. Our preferred parametrization is exactly the one 
  for which $h_0=0$; see Figure \ref{fig:bifurcationstructure}.
 (Also, every component which bifurcates from a period $n$ component
  has period $qn$ for some $q\geq 2$, and the bifurcation angle is
  $p/q\in\Q\setminus\Z$.) 

 The problem with this description is that we do not yet know that 
  $\Psi_W(2\pi i h)\in\C$ for all $h$. However, the address 
  $\extaddr(W,h)$ can be defined combinatorially for all 
  $h\in \Q\setminus\Z$, regardless of whether $\Psi_W(2\pi i 
  h)\in\C$  (and can be computed by a simple algorithm 
  \cite[Lemma 6.3]{expcombinatorics}); 
  this yields the following fundamental 
  result. 

\begin{prop}[Bifurcation Structure {\cite[Theorem 6.8]{expcombinatorics}}] 
   \label{prop:bifurcationstructure}
 Let $W$ be a hyperbolic component and $\s:=\extaddr(W)$. 
  If $n\geq 2$, then there exists a
 unique choice of
 $\Psi_W$ --- called the \emph{preferred
 parametrization} of $W$ --- and a unique map
    $\extaddr(W,\cdot):\Q\setminus\Z\to\Sequ$
 with the following properties.
  \begin{enumerate} 
   \item $\extaddr(W,\cdot)$ is strictly increasing on
         $\{h>0\}$ and (separately) on
         $\{h<0\}$;
   \item $\extaddr(W,\frac{p}{q})$ is an intermediate external address of
           length $qn$; \label{item:intermediateaddress}
   \item if $h\in\Q\setminus\Z$ is such that
          $\Psi_W(2\pi i h)\in\C$, then the parameter $\Psi_W(2\pi i h)$ lies on
          the boundary of $\child:= \Hyp{\extaddr(W,h)}$;
   \item $\displaystyle{\cl{\W(W)}=\cl{\bigcup_{h\in\Q\setminus\Z}
          \W(\child)}}$.
             \label{item:subwakesfillwake}
   \item $\displaystyle{\lim_{h\to+\infty}
          \extaddr(W,h) = \lim_{h\to-\infty} \extaddr(W,h)= \s}$;
          \label{item:addressofW}
   \item $\displaystyle{\lim_{h\nearrow 0} \extaddr(W,h)=\s^+}$ and
         $\displaystyle{\lim_{h\searrow 0} \extaddr(W,h)=\s^-}$.
         \label{item:charaddresses}
  \end{enumerate}

 If $n=1$, then the preferred parametrization is, by definition, the
 map $\Psi_W(z)=z-\exp(z)$, which
 maps $\R^-$ to $(-\infty,-1)$. There exists a unique map
 $\extaddr(W,\cdot)$ which is strictly increasing on all of
 $\Q\setminus\Z$ and satisfies properties
 (\ref{item:intermediateaddress}) to (\ref{item:subwakesfillwake}) above. 
\end{prop}
\begin{remark} From now on, we shall always fix
  $\Psi_W$ to be the preferred parametrization of $W$ (note that this
  is the
  same preferred choice as described in 
  \cite[Theorem 7.1]{expattracting}).  
  In particular, with this definition
  the internal ray of height $h$ is uniquely defined (recall the
  remark after Definition \ref{defn:internalrays}).
\end{remark}

\begin{defn}[Types of Hyperbolic Components]  \label{def:satellite}
 The components $\child$ are called the \emph{child components}
 of $W$. If $W$ is a child component of some other hyperbolic
 component, it is also called a \emph{satellite component}; otherwise
 we say that $W$ is a \emph{primitive component}. 
\end{defn}
\begin{remark} The term \emph{satellite} is more commonly used for a
  hyperbolic component which shares a parabolic boundary point
  with a hyperbolic component of lesser period. Theorem
  \ref{thm:noncentrallanding} shows
  that these definitions are equivalent; compare also Corollary
  \ref{cor:satellite}.
\end{remark}

\begin{prop}[Child Components
   {\cite[Proposition 6.2]{expcombinatorics}}] \label{prop:childcomponents}
 Let $W$ be a hyperbolic component, and let $V:=\child$ be a child component
 of $W$. Then
    $\Psi_W(2\pi ih) = \Psi_{V}(0)$.
 If $\kappa_0:=\Psi_W(2\pi ih)\in\C$, 
 then $W$ and $V$ are the only
  hyperbolic components containing $\kappa_0$ on their boundary. 

 Furthermore, if two components have a common
 parabolic boundary point, then one is a child component of the
 other; no two hyperbolic components have a common child
 component. In particular, no two hyperbolic components of equal
 period have a common parabolic boundary point. \qedd
\end{prop}
\begin{remark} We will show in Proposition \ref{prop:intersection} that no
 two components ever touch at irrational boundary points.
\end{remark}

If $W$ is a hyperbolic component, and $V$ is a child component of $W$,
 then $\W(V)$ is called a \emph{subwake} of $W$. By Proposition
 \ref{prop:bifurcationstructure} (\ref{item:subwakesfillwake}),
 these subwakes exhaust most of the wake of $W$. The following
 statement makes
 this more precise.

\begin{prop}[Subwakes Fill Wake {\cite[Corollary 6.9]{expcombinatorics}}]
  \label{prop:subwakesfillwake}
 Let $W$ be a hyperbolic component of period $n$, and let
 $\s\in\W(W)\setminus\{\extaddr(W)\}$. If $\s$ is not contained in
 any subwake of $W$, then $\s$ is a bounded infinite external address.\qedd
\end{prop}

 We will also need a
  bound on parameters in a given wake, which is given by the following result.

\begin{lem}[Bound on Parameter Wakes
             \protect{\cite[Corollary 4.8]{topescapingnew}}]
  \label{lem:parameterraybound}
 Let $\kappa\in\C$, and suppose that two dynamic rays $g_{\s^1}$ and
 $g_{\s^2}$ have a common landing point. Suppose that there are $n\in\N$
 and
 $M > F^{n-1}(6)$ such that, for every $k> 0$ and $j\in\{1,2\}$,
   $\max_{k\leq \ell < k+n} 2\pi |s^j_{\ell}| \geq M$.
 Then $|\kappa| > F^{-(n-1)}(M-\pi)-2$.
\end{lem}
\begin{proof} We claim first that there is an $\ell\in\{1,2\}$ 
   such that the ray $g_{\s^{\ell}}$
   contains a point which maps into the half strip
   $\{z\in\C:\re z \leq \re\kappa\text{ and }
              |im z|\leq |\im\kappa|+2\}$ under iteration. Indeed, since
   $\s^1\neq\s^2$, there exists some $j\geq 0$ such that the first
   entries of $\sigma^j(\s^1)$ and $\sigma^j(\s^2)$ differ. In order
   for the two rays $g_{\sigma^j}(\s^1)$ and $g_{\sigma^j}(\s^2)$ 
   to land together, at least one of them must
   cross, or land on, one of the lines $\{\im z =
   (2k-1)\pi\}$. Since each of these lines maps to
   $\{\kappa+t:t\in\R^-\}$, the claim follows.

  Let $t_0\geq t_{\s^{\ell}}$ be maximal such that there
  is an $m\geq 0$ with 
  $ \re \Ek^m(g_{\s^{\ell}}(t_0))) \leq \re \kappa$ and
  $ |\im z|\leq |\im \kappa|+2$. 
  Then, for every $j\geq 0$, the forward images of the piece
  $g_{\sigma^j(\s^{\ell})}\bigl([F^j(t_0),\infty)\bigr)$ do not intersect
  any strip boundaries; thus this piece is contained in the strip
   \[ \{z\in\C:\im z \in
                   \bigl((2s_{j+1}-1)\pi,(2s_{j+1}+1)\pi\bigr)\}. \] 

  Let $z := \Ek^m(g_{\s^{\ell}}(t_0))$. By choice of $z$ and the first
  part of Lemma \ref{lem:expogrowth}, 
    $|\im \Ek^j(z)| \leq F^{j}(|\kappa|+2)$
  for all $j$. On the other hand, by assumption there exists some
   $j\leq n-1$ with $2\pi |s^{\ell}_{m+j+1}| \geq M$, and thus 
     \[  M - \pi \leq |\im \Ek^j(z)| \leq F^j(|\kappa|+2)\leq
        F^{n-1}(|\kappa|+2). \qedhere \]
 \end{proof}


\section{Cutting Inside Wakes} \label{sec:squeezingc}

We begin this section by showing that, for every hyperbolic component,
 every internal ray lands in $\C$, except possibly the ray at height
 $h=0$. This suffices to prove our main result (Theorem
 \ref{mainthm:boundary}) for all satellite components
 (see Definition \ref{def:satellite}). We conclude by proving
 part (c) of the Squeezing Lemma; i.e.\ the case of 
 intermediate external addresses.

 The main technique applied in this section
  is to use bifurcating components to cut a given curve
  off from $\infty$. Thus we need to
  make sure that sufficiently many such
  components exist. Luckily this is the case:

\begin{prop}[Bifurcation Angles are Dense] \label{prop:densebifurcations}
 Let $W$ be a hyperbolic component, and let $\Psi_W:\H\to W$ be its
 preferred parametrization. Then the set
  $\{h\in \Q\setminus\Z : \Psi_W(2\pi i h)\in\C \}$
 is dense in $\R$.
\end{prop}
\begin{proof} The set 
   $\{ h \in \R: \Psi_W(2\pi i h)\in\C\}$
  is open by continuity of $\Psi_W$ (Lemma
  \ref{lem:localconnectivity}), and dense by the F.~and M.~Riesz
  theorem. \end{proof}

We can now formulate the main tool we will use in this and the following
 section to transfer combinatorial information to the parameter plane. 
 (Compare Figure \ref{fig:cuttinginwake}.)
\begin{prop}[Curves in the Wake of a Hyperbolic Component]
  \label{prop:cuttinginwake}
 Let $\gamma:[0,\infty)\to\C$ be a curve with $\gamma(t)\to\infty$
 as $t\to\infty$, and suppose that $\gamma$ contains no
 indifferent parameters. If $W$ is a hyperbolic component with
 $\extaddr(\gamma)\in \W(W)$, then there exists a curve in $\C$ which
 separates $\gamma$ from all hyperbolic components $\Hyp{\t}$ with
 $\t\notin \W(W)$. This curve can be chosen to consist
 only of parameters in $W$, parameters in a single 
 child component of $W$ and one
 parabolic parameter on $\partial W$.
\end{prop}
\begin{proof}
 To fix ideas, let us assume that $\gamma$ tends
 to $\infty$ either in or below the preferred homotopy class of
 $W$. By assumption, $\extaddr(\gamma) > \s^-$, so by 
 Proposition \ref{prop:bifurcationstructure}
 (\ref{item:charaddresses}), there exists an $\eps>0$ such that
 $\extaddr(W,h)<\extaddr(\gamma)$ for all rational $h\in(0,\eps)$.
 By Proposition \ref{prop:densebifurcations}, there exists some
 rational $h\in(0,\eps)$
 with $\kappa_0:= \Psi_W(2\pi i h)\in\C$.
 Since $\gamma$ intersects at most one
 hyperbolic component, we can also assume that $\gamma$ does not
 intersect the child component
 $V:=\child$. The desired curve is given by
   \[ \IR \cup \{\kappa_0\} \cup \IRWH{V}{0}. \qedhere \]
\end{proof}

\begin{thm}[Landing of Non-Central Internal Rays] \label{thm:noncentrallanding}
 Let $W$ be a hyperbolic component and let $\Psi_W$ be its preferred
 parametrization. Then $\Psi_W(2\pi i h)\in\C$ for all
 $h\neq 0$.
\end{thm}
\begin{proof}
  Suppose that $h\neq 0$ and $\Psi_W(2\pi i
  h)=\infty$. Let us assume that
   $h > 0$; the case $h<0$ is completely analogous. Consider the internal ray
   \[ \IR:(-\infty,0)\to\C, t\mapsto \Psi_W( t+2\pi i h). \]
 Then the curve $\IR:(-\infty,-1]\to\C$ tends to $\infty$ in the
  preferred homotopy class of $W$ as $t\to-\infty$
  and thus has external address
   $\t:=\extaddr(W)$. Since $\lim_{t\to 0}\IR(t)=+\infty$ by assumption,
   the curve $\IR:[-1,0)\to\C$ also defines an external address
     \[ \s := \extaddr\Bigl(\IR\bigl([-1,0)\bigr)\Bigr). \] 

 By Proposition \ref{prop:densebifurcations}, there exists a rational
  $h_0>h$ for which $\kappa_0 := \Psi_W(t+2\pi i
  h_0)\in\C$. Since $\IR$ surrounds $\kappa_0$,
  it also surrounds the child component $\childH{h_0}$. This shows that $\s < \t$.

 Similarly, there exists some
  $h_1\in\Q\setminus\Z$ between $0$ and $h$ with $\kappa_1:=\Psi_W(2\pi i h_1)
  \in \C$. As in the proof of Proposition \ref{prop:cuttinginwake},
  the curve $\IRH{h_1}\cup \{\kappa_1\} \cup
             \IRWH{\childH{h_1}}{0}$
  surrounds $\IR$. It
  follows that
   $\extaddr(\childH{h_1}) \leq \s \leq \extaddr(W)$; 
   in particular $\s\in \W(W)$.

 By part (b) of the Squeezing Lemma (Theorem \ref{thm:squeezingb}),
  $\s$ is either intermediate or infinite and not exponentially
  bounded; in particular $\K(\s)$ is not periodic. By Corollary
  \ref{prop:subwakesfillwake},  
  there exists some
  child component $V$ of $W$ such that
  $\s\in\W(V)$. By Proposition \ref{prop:cuttinginwake}, there is a
  curve which is disjoint from $W$ and separates the piece
  $\IR([-1,0))$ from $W$. This is a contradiction. \end{proof}

\begin{cor}[Boundaries of Satellite Components]  \label{cor:satelliteboundary}
 Suppose that $V$ is a satellite component.
 Then $\Psi_V:\cl{\H}\to\cl{W}$ is a homeomorphism, and $\partial
 V\cap\C$ is a Jordan arc.
\end{cor}
\begin{proof}
  By Theorem \ref{thm:noncentrallanding}, it only remains to show
 that $\Psi_V(0)\in\C$. Since $V$ is a satellite component,
 there exists a hyperbolic component $W$ and
 $h\in\Q\setminus\Z$
 such that $V=\child$. By Proposition
 \ref{prop:childcomponents} and Theorem \ref{thm:noncentrallanding},
 $\Psi_V(0)=\Psi_W(2\pi ih)\in\C$. \end{proof}

We will now prove part (c) of the Squeezing Lemma.

\begin{thm}[Squeezing Lemma, Part (c)] \label{thm:squeezingc}
 Suppose that $\gamma:[0,\infty)\to\C$ is a curve to $\infty$ which
 does not contain indifferent parameters. If $\s := \extaddr(\gamma)$ is
 intermediate, then $\gamma$ lies in the preferred homotopy class of
 $\Hyp{\s}$.
\end{thm}
\begin{proof} By the discussion preceeding Theorem \ref{thm:squeezinglemma}, 
 the case $\s=\infty$ is trivial. 
 Let us then suppose, by contradiction,
 that $\s\neq\infty$ and that
 $\gamma$ does not lie in the preferred homotopy class of
 $W:=\Hyp{\s}$. To fix our ideas, let us assume that $\gamma$ tends to
 $\infty$ below this homotopy class.

 If $\gamma\subset W$, then by
 Theorem \ref{thm:noncentrallanding},
 $\Phi_W(\gamma(t))\to 0$ as $t\to\infty$, and for any $h\in
 \Q\setminus\Z$ with $h>0$ it follows that
 $\child$ tends to $\infty$ between the preferred
 homotopy class of $W$ and $\gamma$, a contradiction.

 So let us assume that
 $\gamma$ does not intersect $W$.
 For every $j\in\N$, set $\kappa_j := \Psi_W(2\pi i
 (j+\onehalf))\in\C$. 
 Denote by $\Gamma^1_j$ the internal ray of $W$ landing at
 $\kappa_j$ and by $\Gamma^2_j$ the central internal ray of the child
 component $V_j$ bifurcating from $\kappa_j$. Then the curve
  \[ \Gamma_j := \Gamma^1_j \cup \{\kappa_j\} \cup \Gamma^2_j \]
 surrounds $\gamma$. 
 We will derive a contradiction by showing that
  \[ \lim_{j\to\infty} \min \{|\kappa|: \kappa\in \Gamma_j\}= \infty. \]

 Since $\Psi_W$ is continuous in $\infty$, it is clear that
  \[ \lim_{j\to\infty} \inf \{|\kappa|:\kappa\in \Gamma_j^1\} = \infty. \]
 It is thus sufficient to concentrate on the curves $\Gamma_j^2$. Denote
  the characteristic addresses of $V_j$ by $\r_j$ and $\addrt_j$. 
  These addresses are periodic of period $2n$ 
  (where $n$ is the period of $W$),
  and the corresponding dynamic rays land together for every
  parameter on $\Gamma_j^2$. Since all $\r_j$ and $\addrt_j$ are 
  different,
  their largest entries must tend to $\infty$ as $j$ becomes large.
  By Lemma \ref{lem:parameterraybound}, this
  implies that 
  $\inf \{|\kappa|:\kappa\in \Gamma_j^2\}$
  also tends to $\infty$. \end{proof}


\section{Pushing Wakes to Infinity}
 \label{sec:squeezinga}

 In this section, we complete the proof of the Squeezing Lemma by
 showing that the address of any curve $\gamma$ 
 in exponential parameter space
 must be intermediate or exponentially bounded.  In order to do this,
 we shall surround $\gamma$ by curves in hyperbolic components,
 similarly to the previous section, such that these curves must lie
 farther and farther to the right, which yields a
 contradiction. However, we shall need a theorem relating a given
 external address to the wakes which contain it. The theory of
 \emph{internal addresses}, introduced for polynomials in
 \cite{intaddrnew}, gives exactly such a description. We will
 use the following result (compare \cite[Corollary 1.7]{intaddrnew}).

 \begin{prop}[Finding Hyperbolic Components {\cite[Corollary 7.11]{expcombinatorics}}] 
    \label{prop:internaladdress}
   Let $\s$ be an unbounded infinite external address, 
   and suppose that $W$ is a hyperbolic component with
   $\s\in\W(W)$. Let $k$ be maximal such that the first
   $k-1$ entries of $\addu:=\K(\s)$
   and $\KS(W)$ coincide (with the convention that $k=1$ if
   $W=\Hyp{\infty}$ is the period $1$ component). If $\u_k\in\Z$, then
   there exists a hyperbolic component $V$ with $\s\in\W(V)$ and
   $\KS(V)=\per{\u_1 \dots \u_k}$. \qedd
 \end{prop}

 \begin{thm}[Squeezing Lemma, Part (a)]
   \label{thm:squeezinga}
  Suppose that $\gamma:[0,\infty)\to\C$ is a curve to $\infty$ which
   does not contain indifferent parameters. Then $\extaddr(\gamma)$ is
   either intermediate or exponentially bounded.
 \end{thm}
 \begin{proof} Suppose, by contradiction, that 
   $\s:=\extaddr(\gamma)$ is infinite
   but not exponentially bounded. The first step in our proof
   is to find a sequence $W_k$ of hyperbolic components with
   $\s\in\W(W_k)$ such that the combinatorics of $W_k$ becomes
   ``large'' in a suitable sense.

  Set $\addu:=\K(\s)$; then the sequence
   $\addu$ is also
   not exponentially bounded. 
   We can thus choose a subsequence $(\u_{n_k})$ of entries such that
   $|\u_j| < |\u_{n_k}|$ whenever $j<n_k$ and such that
   \[ F^{-n_k}(2\pi|\u_{n_k}|)\to \infty. \]
   We claim that, for every $k$, there exists a hyperbolic component 
    $W_k$ with
   $\s\in\W(W_k)$ and
   \[ \KS(W_k) = \per{\u_1\dots \u_{n_k}}. \]

   Indeed, let $m<n_k$ be maximal such that there exists a hyperbolic
   component $U$ with $\s\in\W(U)$ with $\KS(U)=\per{\u_1\dots \u_m}$. 
   (The existence of such an $m$ 
    follows from Proposition \ref{prop:internaladdress}.) 
   By
   maximality of $m$ and Proposition \ref{prop:internaladdress},
   $\KS(U)$ and $\addu$ agree on the first $n_k-1$ entries. On the other
   hand, all entries of $\KS(U)$ are different from $\u_{n_k}$. The
   existence of a component $W_k$ with the desired properties now
   follows by applying Proposition \ref{prop:internaladdress} to $\s$
   and $U$.

 By Proposition \ref{prop:cuttinginwake}, there is a
  curve $\Gamma_k\subset\C$ which tends to $+\infty$ in both directions,
  which separates $\gamma$ from a left half plane, and which
  consists only of parameters in $W_k$ and a child component of $W_k$
  together with a common boundary point. We will show that the
  curves $\Gamma_k$ tend to $\infty$ uniformly as $k\to\infty$. 

 Let $\r^k$ and $\addrt^k$ denote the characteristic addresses
  of $W_k$. Since $\KS(W_k)=\per{\u_1\dots \u_{n_k}}$, for every $j$ 
  the $jn_k$-th
  entries of $\r^k$ and $\addrt^k$ are of size at least
  $|\u_{n_k}|-1$. By Lemma \ref{lem:parameterraybound}, it follows
  that
  \[ \Gamma_k \subset \bigl\{\kappa: |\kappa| >
  \F^{-n_k+1}(2\pi|\u_{n_k}|-3\pi)-2\bigr\}. \]
  On the other hand,
  $F^{-n_k+1}(2\pi|\u_{n_k}|-3\pi)\to\infty$. This
  is a contradiction because every $\Gamma_k$ surrounds $\gamma$.
  \end{proof}

This completes the proof of the Squeezing Lemma, and thus of Theorem
\ref{mainthm:boundary}. \qedoutsideproof

\section{Further Results} \label{sec:furtherresults}

With Theorem \ref{mainthm:boundary} now proved, we 
discuss a number of important properties of the bifurcation structure 
in exponential parameter space: we prove that boundaries of 
hyperbolic components are analytic except at roots and co-roots, we 
solve another conjecture of Eremenko and Lyubich on "bifurcation 
trees", and we give an intrinsic definition of wakes in parameter 
space using landing properties of dynamic and parameter rays. 

\subsection*{Analyticity of boundaries}
Let $W$ be a hyperbolic component of period $n$. We will call
$\Psi_W(0)$ the \emph{root} of $W$; the points of $\Psi_W(2\pi i
(\Z\setminus\{0\})$ are called \emph{co-roots}. 
We will show that the boundary of $W$ is an
analytic curve, except at its co-roots and possibly at its root point. Note
that it is a priori clear that $\partial W$ is a piecewise analytic
curve. We need to rule out the existence of critical points
of the multiplier map $\mu$ on $\partial W$.

\begin{prop}[Closures Intersect at Parabolic Points]
 \label{prop:intersection}
 Let $W$ be a hyperbolic component and let $\kappa\in \partial W$ be
 an irrationally indifferent parameter. Then no other hyperbolic component contains
 $\kappa$ on its boundary. 
\end{prop}
\begin{proof} Suppose that $W_1$ and $W_2$ are two hyperbolic components
 which have an irrational boundary parameter $\kappa_0$ in
 common. Let $h_0\in\R\setminus\Q$ with
 $\kappa_0=\Psi_{W_1}(ih_0)$, and suppose, to fix our ideas, that
 $h_0>0$. Then for every $h\in (0,h_0)\cap\Q$, the curve
  \[ \IRW{W_1} \cup \{\Psi_{W_1}(ih)\}\cup
     \IRWH{\childW{W_1}}{0} \]
 separates $\kappa_0$, and thus also
  $W_2$, from every component which does not lie in the wake
  of $W_1$ (compare Proposition \ref{prop:cuttinginwake}). 

 This proves
  that $\extaddr(W_2)\in\W(W_1)$. 
  By symmetry, also $\extaddr(W_2)\in
  \W(W_1)$, which is a contradiction. \end{proof}

\begin{cor}[Analytic Boundary] \label{cor:analyticboundary}
 Let $W$ be a hyperbolic component.
 Then the  function $h\mapsto \Psi_W(2\pi i h)$ is analytic in
 $\R\setminus\Z$. Furthermore, $\partial W$ has a
 cusp in every co-root of $W$. The boundary is analytic or has a cusp
 in the root of $W$ depending on whether $W$ is a satellite or
 primitive component, respectively.
\end{cor}
\begin{proof}
 Let $h\in\R\setminus\Z$; then the multiplier map $\mu$ extends analytically
 to a neighborhood of $\kappa:=\Psi_W(2\pi i h)$. Since
 $\mu\circ\Psi_W=\exp$, analyticity of $\Psi_W$ follows
 unless $\mu$ has a critical point in $\kappa$. 
 In that case, let $D$ be a small disk around $e^{2\pi i h}$. The
 preimage of $D\cap\D$ under $\mu$ has at least two components
 $D_1,D_2$ whose boundary contains
 $\kappa$. By Proposition
 \ref{prop:childcomponents} (if $\theta\in\Q$) resp.~Proposition
 \ref{prop:intersection} (otherwise), $W$ is the only component of
 period at most $n$ containing $\kappa$ in its boundary, and thus
 $D_1\cup D_2\subset W$. However, this is impossible by Corollary
 \ref{cor:landingpoints}. 

The statement about points of $\Psi_W(2\pi i \Z)$ is proved in a
 similar way as for the Mandelbrot set \cite[Lemmas 6.1 and
 6.2]{jackrays}. By an elementary local argument, for 
 every root or co-root $\kappa_0$ there is a neighborhood $U$ of $\kappa_0$ 
 such that the multiplier map is defined at least on some double cover 
 of $U\setminus\{\kappa_0\}$ (using the fact that every primitive parabolic 
 parameter is at most a double parabolic, because there is only one 
 singular orbit). In the satellite case, the lower-period multiplier 
 is defined in a neighborhood of $\kappa_0$ in the $\kappa$-plane and 
 has no critical point at $\kappa_0$, so $W$ has analytic boundary. It 
 follows that the higher-period multiplier has similar properties. Otherwise,
 $\kappa_0$ is on the boundary of a 
 single hyperbolic component, so $\mu$ must be injective in a 
 neighborhood of $\kappa_0$ on the double cover of $U$, and the claim 
 follows as in \cite{jackrays}.
\end{proof}

\subsection*{Bifurcation trees} 
A consequence of
 Proposition \ref{prop:intersection} is the following
 characterization of satellite hyperbolic components.
\begin{cor}[Satellite Components] \label{cor:satellite}
 Let $W$ be a hyperbolic component. Then the following are
  equivalent
  \begin{enumerate}
   \item $W$ is a satellite component.
   \item There is a hyperbolic component $V\neq W$ with
     $\Psi_W(0)\in \partial V$.
   \item There is a hyperbolic component $V$ of period less than
     $W$ such that $\partial V\cap\partial W\cap\C\neq\emptyset$.
  \end{enumerate}
\end{cor}
\begin{proof} By Proposition \ref{prop:childcomponents},
  the root of every satellite
  component lies on a hyperbolic component of smaller period. 
  If $\Psi_W(0)$ lies on the boundary of another
  hyperbolic component $V$, then it follows from Proposition
  \ref{prop:childcomponents} that $W$ is a child component of $V$, 
  and thus the period of $V$ is smaller than that of $W$.
  Finally, suppose that $V$ is a component of smaller period than $W$ 
  such that
  $W$ and $V$ have a common finite boundary point. Then
  this boundary point is parabolic by Proposition
  \ref{prop:intersection}. Thus, by Proposition \ref{prop:childcomponents},
  $W$ is a child component of $V$. \end{proof}

The \emph{bifurcation forest} of hyperbolic components is the
(infinite) graph with one vertex for each
hyperbolic component and one edge for each pair of components whose
closures have a finite intersection point. A \emph{bifurcation
tree} is any component of this graph. 
The periods of hyperbolic components make every bifurcation tree an oriented
tree, which thus has a unique root point of lowest period. By Corollary
 \ref{cor:satellite}, all vertices of 
the tree correspond to satellite components, with the exception of 
the root point of the tree, which always is a primitive component. 
Conversely, every primitive component is the root point of its own 
bifurcation tree. Moreover, different primitive components have disjoint 
bifurcation trees.
It was conjectured in
\cite{alexmisharussian} that there are infinitely many bifurcation
trees. We will now prove this fact.

\begin{cor}[Infinitely Many Bifurcation Trees]  
  \label{cor:infinitelymanybifurcationtrees} 
 There are infinitely many bifurcation trees.
\end{cor}
\begin{proof}
  It suffices to prove that there are infinitely many primitive components.
  We will show that, for every $k>0$,
  the component
  $\Hyp{0(k+\onehalf)\infty}$ is primitive. Indeed, note that
  $0(k+\onehalf)\infty\in \W(\Hyp{\onehalf\infty})=(\per{01},\per{10})$.
  Thus $\Hyp{0+\onehalf)\infty}$
  does not bifurcate from the unique period $1$ component
  $\Hyp{\infty}$,
  and is thus primitive.
  \end{proof}

\subsection*{Wakes and periodic parameter rays}
 In \cite[Corollaries IV.4.4 and IV.5.2]{habil}, 
 it was shown (without using the results of this
 article), that every parabolic parameter is the landing point of
 either one or two parameter rays at periodic addresses. (Recall the
 definition of the parameter rays $\PR$ in Definition \ref{defn:parrays}.)
 More precisely, suppose that
 $\kappa$ is a parabolic parameter; say
 $\kappa=\Psi_W(2\pi i k)$ for some hyperbolic component $W$ of period
 $n$ and some $k\in\Z$. If $k=0$, then $\kappa$ is the landing point
 of exactly two periodic parameter rays, namely those at the characteristic
 addresses of $W$, which both have period $n$. 
 If $k\in\Z\setminus\{0\}$, then $\kappa$ is the
 landing point of a single periodic parameter ray, at
 the address
   $\extaddr(W,k) := \lim_{h\to k} \extaddr(W,h)$. These
 addresses are called \emph{sector boundaries} of $W$.

 Using Theorem \ref{mainthm:boundary} and the fact that every
 periodic external address is a characteristic address or a
 sector boundary \cite[Lemma 7.4]{expcombinatorics}, 
 we obtain the following result \cite[Theorem V.7.2]{habil}.
 \begin{thm}[Periodic Parameter Rays Land]
  \label{thm:periodicparameterrays}
  Every periodic parameter ray lands at a parabolic
 parameter. Conversely, every parabolic parameter is the landing point
 of either one or two periodic parameter rays. \qedd
 \end{thm}

 We have so
 far only defined the \emph{combinatorial}
 wake of a hyperbolic component $W=\Hyp{\s}$,
 as the interval $(\s^-,\s^+)$. 
  However, Theorem \ref{thm:periodicparameterrays} (and the
 discussion preceding it) suggest the following definition of the wake
 as a subset of parameter space (which is
 analogous to the usual definition of wakes in the Mandelbrot set).
 \begin{defn}[Wake] \label{defn:parameterwake}
  Let $W$ be a hyperbolic component with characteristic addresses
   $\s^-$ and $\s^+$. Then the \emph{(parameter) wake} 
  of $W$ is defined to be the
  component of 
   $\C\setminus \bigl( \PRS_{\s^-}\cup \PRS_{\s^+} \cup
      \Psi_W(0)\bigr)$
  which contains $W$.
 \end{defn}
 For the Mandelbrot set, the parameter
 wake of a hyperbolic component $W$ coincides
 with the set of all parameters for which the dynamic rays at the
 characteristic addresses of $W$ land at a common repelling periodic
 point. For exponential maps, the equivalence of this definition to
 the previous one ---
 and the landing of periodic dynamic rays for all expoential maps
 with nonescaping singular value ---
 was recently shown in \cite{landing2new}, 
 using Theorem \ref{thm:periodicparameterrays} and holomorphic
 motions.
 \begin{thm}[Characterization of Wakes]
   Let $W$ be a hyperbolic component with characteristic addresses
   $\s^-$ and $\s^+$. Then the parameter wake of $W$
  coincides with the set of all
  parameters $\kappa$ for which $g_{\s^-}$ and $g_{\s^+}$ have a common
  repelling landing point. \qedd
 \end{thm}

\subsection*{Nonhyperbolic components}
The main open question about exponential parameter space, like for 
quadratic polynomials and the Mandelbrot set, is whether hyperbolic 
dynamics is dense. In other words, we need to show that there are no 
non-hyperbolic (or "queer") stable 
components. It seems possible to show at 
least that any non-hyperbolic component must be bounded: the first
combinatorial step is to show that hyperbolic components, as well as 
parameter rays at periodic and preperiodic external addresses 
together with their landing points, disconnect the $\kappa$-plane 
into complementary components such that each component 
is separated from all external addresses with at most two exceptions, 
both of which must be exponentially bounded. In particular, if a 
non-hyperbolic component is unbounded, it must ``squeeze'' to $\infty$ 
very close to one or two parameter rays, and a stronger variant of 
the Squeezing Lemma might prevent this from happening. This would 
show that all non-hyperbolic components were bounded. 
However, most features in 
exponential parameter spaces are unbounded (for example, all 
hyperbolic components and all puzzle pieces). If one could prove that 
every non-hyperbolic component had to be unbounded as well, this 
would settle density of hyperbolicity in an interesting way.


 \nocite{krantz}
\bibliographystyle{hamsplain}
\small{\bibliography{/Latex/Biblio/biblio}}

\end{document}

%% file: symboldefs.tex
\renewcommand{\H}{\mathbb{H}}

\renewcommand{\theta}{\vartheta}
\renewcommand{\phi}{\varphi}
\renewcommand{\rho}{\varrho}

\newcommand{\Ek}{E_{\kappa}}

\newcommand{\onehalf}{\frac{1}{2}}

\newcommand{\ul}[1]{\underline{#1}}

\newcommand{\Sequ}{\mathcal{S}}

\newcommand{\Sequb}{\overline{\Sequ}}

\newcommand{\adds}{\underline{s}}

\newcommand{\s}{\adds}

\newcommand{\F}{\mathcal{F}}

\newcommand{\ts}{t_{\s}}

\renewcommand{\u}{{\tt u}}

\newcommand{\bdyit}[2]
             {{\rule{0pt}{0pt}_{\mbox{$\scriptstyle #2$}}^{\mbox{%
                   $\scriptstyle #1$}} }}
\renewcommand{\j}{{\tt j}}
\newcommand{\itj}{\bdyit{\j}{\j-1}}

\newcommand{\rt}{\tilde{r}}
\newcommand{\addrt}{\ul{\tilde{r}}}
                               
\newcommand{\addu}{\ul{\u}}

\newcommand{\addt}{\ul{t}}
\renewcommand{\t}{\addt}

\renewcommand{\r}{\ul{r}}
\renewcommand{\rt}{\widetilde{\r}}

\newcommand{\extaddr}{\operatorname{addr}}

\newcommand{\itin}{\operatorname{itin}}
\newcommand{\K}{\mathbb{K}}

\newcommand{\KS}{\K^{*}}

\newcommand{\gs}{g_{\adds}}

\newcommand{\periodic}[1]{\overline{#1}}
\newcommand{\per}[1]{\periodic{#1}}

\newcommand{\W}{\mathcal{W}}

\newcommand{\Hplane}{\mathbb{H}^+}

  \newcommand{\B}{\mathcal{B}}

\newcommand{\Hyp}[1]{\operatorname{Hyp}(#1)}

\newcommand{\Bif}{\operatorname{Bif}}
\newcommand{\child}{\Bif(W,h)}
\newcommand{\childH}[1]{\Bif(W,#1)}
\newcommand{\childW}[1]{\Bif(#1,h)}

%% file: bifurcations.bbl
\providecommand{\href}[2]{#2}\def\polhk#1{\setbox0=\hbox{#1}{\ooalign{\hidewid%
th \lower1.5ex\hbox{`}\hidewidth\crcr\unhbox0}}}
  \def\polhk#1{\setbox0=\hbox{#1}{ \ooalign{\hidewidth
  \lower1.5ex\hbox{`}\hidewidth\crcr\unhbox0}}} \input{cyracc.def} \def\j{{\u
  i}} \def\J{{\u I}} \newfont{\cyrit}{wncyi10 at 12pt}\def\cprime{$'$}
\providecommand{\bysame}{\leavevmode\hbox to3em{\hrulefill}\thinspace}
\begin{thebibliography}{BDG2}

\bibitem[AL]{avilalyubichfeigenbaum2}
Artur Avila and Mikhail Lyubich, \emph{{H}ausdorff dimension and conformal
  measures of {F}eigenbaum {J}ulia sets}, Preprint \#2005/05, Institute for
  Mathematical Sciences, SUNY Stony Brook, 2005,
  \mbox{\href{http://www.arXiv.org/abs/math.DS/0408290}{arXiv:math.DS/0408290}%
}.

\bibitem[BR]{bakerexp}
I.~Noel Baker and Philip~J. Rippon, \emph{Iteration of exponential functions},
  Ann. Acad. Sci. Fenn. Ser. A I Math. \textbf{9} (1984), 49--77.

\bibitem[BD]{tying}
Ranjit Bhattacharjee and Robert~L. Devaney,
  \emph{\href{http://math.bu.edu/people/bob/papers/ranjit.ps}{Tying hairs for
  structurally stable exponentials}}, Ergodic Theory Dynam. Systems \textbf{20}
  (2000), no.~6, 1603--1617.

\bibitem[BDG1]{dghnew1}
Clara Bodel{\'o}n, Robert~L. Devaney, Michael Hayes, Gareth Roberts, Lisa~R.
  Goldberg, and John~H. Hubbard,
  \emph{\href{http://math.bu.edu/people/bob/papers/hairs.ps}{Hairs for the
  complex exponential family}}, Internat. J. Bifur. Chaos Appl. Sci. Engrg.
  \textbf{9} (1999), no.~8, 1517--1534.

\bibitem[BDG2]{dghnew2}
\bysame, \emph{\href{http://math.bu.edu/people/bob/papers/hairs-2.ps}{Dynamical
  convergence of polynomials to the exponential}}, J. Differ. Equations Appl.
  \textbf{6} (2000), no.~3, 275--307.

\bibitem[D]{devaneybifurcation}
Robert~L. Devaney, \emph{Julia sets and bifurcation diagrams for exponential
  maps}, Bull. Amer. Math. Soc. (N.S.) \textbf{11} (1984), no.~1, 167--171.

\bibitem[DFJ]{dfj}
Robert~L. Devaney, N{\'u}ria Fagella, and Xavier Jarque,
  \emph{\href{http://math.bu.edu/people/bob/papers/kneadings.ps}{Hyperbolic
  components of the complex exponential family}}, Fund. Math. \textbf{174}
  (2002), no.~3, 193--215.

\bibitem[DGH]{dgh}
Robert~L. Devaney, Lisa~R. Goldberg, and John~H. Hubbard, \emph{A dynamical
  approximation to the exponential map by polynomials}, Preprint, MSRI
  Berkeley, 1986, published as \cite{dghnew1,dghnew2}.

\bibitem[DH]{orsay}
Adrien Douady and John Hubbard, \emph{Etude dynamique des polyn{\^o}mes
  complexes}, Pr{\'e}publications math{\'e}mathiques d'Orsay (1984 / 1985),
  no.~2/4.

\bibitem[E]{alexescaping}
Alexandre~{\`E}. Eremenko, \emph{On the iteration of entire functions},
  Dynamical systems and ergodic theory (Warsaw, 1986), Banach Center Publ.,
  vol.~23, PWN, Warsaw, 1989, pp.~339--345.

\bibitem[EL1]{alexmisharussian}
Alexandre~{\`E}. Eremenko and Mikhail~Yu. Lyubich, \emph{Iterates of entire
  functions}, Preprint, Physico-Technical Institute of Low-Temperature Physics
  Kharkov, 1984; translation in Stony Brook IMS Preprint \#1990/04, 
  published in modified form as \cite{alexmisha}.

\bibitem[EL2]{alexmisha}
\bysame, \emph{Dynamical properties of some classes of entire functions}, Ann.
  Inst. Fourier (Grenoble) \textbf{42} (1992), no.~4, 989--1020.

\bibitem[F]{markus}
Markus F{\"o}rster, \emph{Parameter rays for the exponential family},
  Diplomarbeit, Techn. Univ. M{\"u}nchen, 2003, Available as
  \href{http://www.math.sunysb.edu/cgi-bin/thesis.pl?thesis03-3}{Thesis
  2003-03} on the \href{http://www.math.sunysb.edu/dynamics/theses}{Stony Brook
  Thesis Server}.

\bibitem[FRS]{markuslassedierk}
Markus F\"orster, Lasse Rempe, and Dierk Schleicher, \emph{Classification of
  escaping exponential maps}, Preprint, 2004,
  \mbox{\href{http://www.arXiv.org/abs/math.DS/0311427}{arXiv:math.DS/0311427}%
}, to appear in Proc. Amer. Math. Soc.

\bibitem[FS]{markusdierk}
Markus F\"orster and Dierk Schleicher, \emph{Parameter rays for the exponential
  family}, Preprint, 2005,
  \mbox{\href{http://www.arXiv.org/abs/math.DS/0505097}{arXiv:math.DS/0505097}%
}, to appear in Ergodic Theory Dynam.\ Systems.

\bibitem[GK{\'S}]{nonrecurrentmero}
Jacek Graczyk, Janina Kotus, and Grzegorz {\'S}wi{\polhk{a}}tek,
  \emph{Non-recurrent meromorphic functions}, Fund. Math. \textbf{182} (2004),
  no.~3, 269--281.

\bibitem[K]{krantz}
Steven~G. Krantz, \emph{Function theory of several complex variables}, John
  Wiley \& Sons Inc., New York, 1982.

\bibitem[McM]{mcmullenrenormalization}
Curtis~T. McMullen, \emph{Complex dynamics and renormalization}, Annals of
  Mathematics Studies, vol. 135, Princeton University Press, Princeton, NJ,
  1994.

\bibitem[M1]{jackdynamicsthird}
\bysame, \emph{Dynamics in one complex variable}, third ed., Annals of
  Mathematics Studies, vol. 160, Princeton University Press, Princeton, NJ,
  2006.

\bibitem[M2]{jacktwocriticalpoints}
\bysame, \emph{On rational maps with two critical points}, Experiment. Math.
  \textbf{9} (2000), no.~4, 481--522,
  \mbox{\href{http://www.arXiv.org/abs/math.DS/9709226}{arXiv:math.DS/9709226}%
}.

\bibitem[M3]{jackrays}
\bysame, \emph{Periodic orbits, externals rays and the {M}andelbrot set: an
  expository account}, Ast\'erisque (2000), no.~261, xiii, 277--333,
  \mbox{\href{http://www.arXiv.org/abs/math.DS/9905169}{arXiv:math.DS/9905169}%
}, G\'eom\'etrie complexe et syst\`emes dynamiques (Orsay, 1995).

\bibitem[R1]{landing2new}
Lasse Rempe, \emph{A landing theorem for periodic rays of exponential maps},
  Proc. Amer. Math. Soc \textbf{134} (2006), no.~9, 2639--2648,
  \mbox{\href{http://www.arXiv.org/abs/math.DS/0307371}{arXiv:math.DS/0307371}%
}.

\bibitem[R2]{topescapingnew}
\bysame, \emph{Topological dynamics of exponential maps on their escaping
  sets}, Ergodic Theory Dynam.~Systems \textbf{26} (2006), no.~6, 1939--1975,
  \mbox{\href{http://www.arXiv.org/abs/math.DS/0309107}{arXiv:math.DS/0309107}%
}.

\bibitem[RS]{expcombinatorics}
Lasse Rempe and Dierk Schleicher, \emph{Combinatorics of bifurcations in
  exponential parameter space}, Transcendental Dynamics and Complex Analysis
  (P.~Rippon and G.~Stallard, eds.), London Math. Soc. Lecture Note Ser.,
  Cambridge Univ. Press, 2007,
  \mbox{\href{http://www.arXiv.org/abs/math.DS/0408011}{arXiv:math.DS/0408011}%
}, pp.~317--370.

\bibitem[S1]{habil}
Dierk Schleicher, \emph{On the dynamics of iterated exponential maps},
  Habilitation thesis, TU M\"unchen, 1999.

\bibitem[S2]{expattracting}
\bysame,
  \emph{\href{http://www.math.helsinki.fi/Annales/Vol28/schleich.html}{Attract%
ing dynamics of exponential maps}}, Ann. Acad. Sci. Fenn. Math. \textbf{28}
  (2003), 3--34.

\bibitem[S3]{cras}
\bysame, \emph{Hyperbolic components in exponential parameter space}, C. R.
  Math. Acad. Sci. Paris \textbf{339} (2004), no.~3, 223--228.

\bibitem[S4]{fibers}
\bysame, \emph{On fibers and local connectivity of {M}andelbrot and {M}ultibrot
  sets}, Fractal geometry and applications: a jubilee of Beno\^\i t Mandelbrot.
  Part 1, Proc. Sympos. Pure Math., vol.~72, Amer. Math. Soc., Providence, RI,
  2004, pp.~477--517.

\bibitem[S5]{intaddrnew}
\bysame, \emph{Internal addresses in the {M}andelbrot set and irreducibility of
  polynomials}, Preprint, 2007,
  \mbox{\href{http://www.arXiv.org/abs/math.DS/9411238v2}{arXiv:math.DS/941123%
8v2}}, updated version of Stony Brook IMS Preprint \#1994/19.


\bibitem[SZ1]{expescaping}
Dierk Schleicher and Johannes Zimmer, \emph{Escaping points of exponential
  maps}, J. London Math. Soc. (2) \textbf{67} (2003), no.~2, 380--400.

\bibitem[SZ2]{expper}
\bysame,
  \emph{\href{http://www.math.helsinki.fi/Annales/Vol28/zimmer.html}{Periodic
  points and dynamic rays of exponential maps}}, Ann. Acad. Sci. Fenn. Math.
  \textbf{28} (2003), 327--354.

\bibitem[Sh]{mitsudim}
Mitsuhiro Shishikura, \emph{The boundary of the {M}andelbrot set has
  {H}ausdorff dimension two}, Ast\'erisque (1994), no.~222, 7, 389--405,
  Complex analytic methods in dynamical systems (Rio de Janeiro, 1992).

\bibitem[UZ1]{urbanskizdunik1}
Mariusz Urba{\'n}ski and Anna Zdunik,
  \emph{\href{http://www.math.unt.edu/~urbanski/papers/z2011703.ps}{The finer
  geometry and dynamics of the hyperbolic exponential family}}, Michigan Math.
  J. \textbf{51} (2003), no.~2, 227--250.

\bibitem[UZ2]{urbanskizdunik4}
\bysame, \emph{Instability of exponential {C}ollet - {E}ckmann maps}, Preprint,
  2005, to appear in Israel J. Math.

\bibitem[W]{whyburnanalytictopology}
Gordon~Thomas Whyburn, \emph{Analytic {T}opology}, American Mathematical
  Society Colloquium Publications, v. 28, American Mathematical Society, New
  York, 1942.

\bibitem[Y]{instability}
Zhuan Ye, \emph{Structural instability of exponential functions}, Trans. Amer.
  Math. Soc. \textbf{344} (1994), no.~1, 379--389.

\end{thebibliography}
